\documentclass[12pt]{amsart}

\usepackage{graphicx}
\usepackage{multirow,arydshln}

\usepackage[colorlinks=false, pdfborder={0 0 0}]{hyperref}  
\usepackage{xcolor} 

\newcommand{\ie}{\textit{i.e.}}

\begin{document}

\def\kb{{\mathbf k}}
\def\tb{{\mathbf t}}
\def\Sc{{\mathcal{S}}}
\def\lcm{{\rm{lcm}}}
\def\NN{{\mathbb{N}}}

\newtheorem{thm}{Theorem}[section]
\newtheorem{Proposition}[thm]{Proposition}
\newtheorem{Definition}[thm]{Definition}
\newtheorem{Lemma}[thm]{Lemma}
\newtheorem{Example}[thm]{Example}
\newtheorem{Remark}[thm]{Remark}

\title[A \texttt{C++} class for algebraic reliability computations]{A \texttt{C++} class for multi-state algebraic reliability computations}
\author{A. M. Bigatti}
\address{Universit\'a degli Studi Genova, Italy}
\email{bigatti@dima.unige.it}
\author{P. Pascual-Ortigosa}
\address{Universidad de La Rioja, Spain}
\email{papasco@unirioja.es}
\author{E. S\'aenz-de-Cabez\'on}
\address{Universidad de La Rioja}
\email{eduardo.saenz-de-cabezon@unirioja.es}

\begin{abstract}
  We present the design and implementation of a \texttt{C++} class for reliability analysis of multi-state systems using an algebraic approach based on monomial ideals. The class is implemented
  within the open-source \verb|CoCoALib| library and provides functions to compute system reliability and bounds.
    The algorithms we present may be applied to general systems with independent components having identical or non-identical probability distributions.
\end{abstract}

\maketitle
\section{Introduction}
The development and implementation of efficient algorithms for system reliability computations is an important task in reliability engineering. Many algorithms exist, and are available to the community in a variety of forms. Some are included in large versatile commercial systems \cite{reliasoft,ALD,ITEM}, others are offered as packages, functions or libraries in mathematical software systems of general purpose languages, for example \verb|Matlab| \cite{LC16,R15}, \verb|Python| \cite{R20} or \verb|R| \cite{SYGBL08}. Others still are directly distributed by the authors as stand-alone software, like \texttt{SHARPE} \cite{sharpe}.

In this paper we describe the \verb|C++| class which implements our algebraic approach to system reliability, and can be integrated in other software systems.  In this way it may become available in different forms and toolboxes to researchers, software developers and reliability engineers.
The language \verb|C++| (now in its version 17, standard ISO/IEC 14882) is a widely used \cite{tiobe} object-oriented general purpose computer language, both for very large systems and for small ad-hoc applications. 
Among the virtues of \texttt{C++} is its integrability with other languages and also its high performance, meeting the need for fast and reliable computations.
Our class is implemented in the \verb|CoCoALib| library \cite{CoCoALib}, which is a \texttt{C++} library for Computations in Commutative Algebra, currently at its version \verb|0.99712| (December 2020). It is open source and free.

The main feature of the \verb|C++| class introduced in this paper is that it is applicable to a large variety of systems, with or without  a known identifiable structure, and can be used to compute the reliability (and bounds) of systems having independent identical or non-identical components. The performance of this class is good in terms of time requirements, being able to compute the reliability of systems with hundreds of components and tens of thousands of minimal paths or cuts. Even though there exist optimized algorithms for several kinds of systems which are faster than the ones presented here, ours are useful to analyze systems for which no specialized algorithms are known, and to benchmark new algorithms for particular types of systems.

The outline of the paper is the following. Section \ref{sec:algebraicReliabilityIntro} gives an introduction to the algebraic approach for system reliability which is at the mathematical core of our implementation. We describe the main functions and design of the class in Section~\ref{sec:class}.
Finally, in Section~\ref{sec:examples}, we show some examples of its use and the results of some computer experiments.
All the code of the class and the examples are available at
{\small
  \verb|http://www.dima.unige.it/~bigatti/data| \verb|/AlgebraicReliability/|}.

\section{Algebraic reliability of coherent systems}\label{sec:algebraicReliabilityIntro}

The algebraic approach to system reliability is based on resolutions and Hilbert series of ideals in rings of polynomials in several indeterminates. Although at first
it might  look like an abstract and theoretical method, it is made practically
applicable by the combinatorial nature of monomial ideals and it
is supported by strong results and algorithms 
in the area. Other authors have used algebraic structures in system reliability analysis before, see for instance the seminal works \cite{BP88, BP88b} or the application of algebraic structures to network reliability \cite{S91, HKCD95}.
In particular, the Universal Generating Function method (UGF), introduced in \cite{U87} and described in more detail in books like \cite{L05}, uses the exponents and coefficients of polynomial-like structures to encode the performance distribution and probabilities of multi-state systems. The UGF method is very flexible and has been applied to several types of multi-state systems, see for instance \cite{LSX19} for a recent example.

In this section we give an overview of the basics of the approach by monomial ideals, that is the backbone of the algorithms in the \verb|C++| class we introduce in Section~\ref{sec:class}.

\subsection{Coherent systems}\label{sec:coherentSystems}
A system $\Sc$ consists of $n$ components which are its elementary units, we denote the $n$ system components by $c_i$ with $i\in\{1,\dots,n\}$. At each moment in time the system is in one of a discrete set of levels $S=\{0,1,\dots, M\}$ indicating growing levels of performance or of failure.
Each individual component $c_i$ of the system can be in one of a discrete set of levels $S_i=\{0,\dots,M_i\}$.
A state of a component is its level and a
state of the system is
the $n$-tuple of its components' states. Given two states $s=(s_1,\dots,s_n)$ and $t=(t_1,\dots,t_n)$ we say that $s\geq t$ if $s_i\geq t_i$ for all
$i=1,\dots,n$ and, conversely, that $s\leq t$ if $s_i\leq t_i$ for all $i$. The level of performance of the system is determined in terms of the states of the components by a {\em structure function} $\Phi:S_1\times \cdots \times S_n\longrightarrow S$. The system $\Sc$ is said to be {\em coherent} if $\Phi$ is non-decreasing and each component is relevant to the system, \ie~for each
component $c_i$ there exist a system state $s=(s_1,\dots,s_n)$ 
and two different levels $j,k\in S_i$ such that
$\Phi(s_{i,j})\neq\Phi(s_{i,k})$,
where $s_{i,\ell} = (s_1,\dots,s_{i-1},\ell,s_{i+1},\dots,s_n)$.

We distinguish between {\em working systems}, also called {\em path systems}, \ie~systems described by the working states of their components, and {\em failure systems} or {\em cut systems} \ie~systems described by the failure states of their components\footnote{There are several definitions of multi-state systems in the literature, see the relationships between them in the diagram given by Natvig in \cite{N11} Figure 2.1. The algebraic methodology that we use can be applied to multi-state strongly coherent systems, multi-state coherent systems and multi-state weakly coherent systems as defined in \cite{N11} Definition 2.4 since the algebraic expressions of reliability are not affected by irrelevant components for each of the levels of performance of the system. }.

In {\em path systems} the focus is on working states. In this case the levels of the system, $\{0,\dots,M\}$ indicate growing levels of performance, the system being in level $0$ indicates that the system is failing, and level $j>i$ indicates that the system is performing at level $j$ better than at level $i$.
For  each component $c_i$ and for each of its levels $j$,
we denote by $p_{i,j}$ the probability that  $c_i$ is performing at
level $\ge j$.
The $j$-reliability of $\Sc$, denoted by $R_j(\Sc)$ is the probability that $\Sc$ is performing at level $\ge j$; conversely, the $j$-unreliability of $\Sc$, denoted $U_j(\Sc)$, is $1-R_j(\Sc)$. A {\em path system} is given at level $j$ by its set of $j$-working states \ie~those tuples $(s_1,\dots,s_n) \in S_1\times\cdots\times S_n$ such that $\Phi(s_1,\dots,s_n)\geq j$. We say that a state $(s_1,\dots,s_n)\in S_1\times\cdots\times S_n$ is a {\em minimal $j$-working state} or {\em minimal $j$-path} if  $\Phi(s_1,\dots,s_n)\geq j$ and $\Phi(t_1,\dots,t_n)< j$ whenever all $t_i\leq s_i$ and at least in one case the inequality is strict. We say that a state $(s_1,\dots,s_n)\in S_1\times\cdots\times S_n$ is a {\em minimal $j$-failure state} or {\em minimal $j$-cut} if  $\Phi(s_1,\dots,s_n)< j$ and $\Phi(t_1,\dots,t_n)\geq j$ whenever all $t_i\geq s_i$ and at least one of the inequalities is strict. Path systems are usually denoted by :G (for {\em good}) in the literature.

\begin{Example}\label{ex:seriesG}
  Let $\Sc$ be a path system with $3$ components $c_1,c_2,c_3$ such that $S=S_1=S_2=S_3=\{0,1,2\}$. The structure function of $\Sc$ is given by $\Phi(s_1,s_2,s_3)=\min\{s_1,s_2,s_3\}$ \ie~$\Sc$ is a multi-state {\em series:G} system, the system works at level $j>0$ only if all of its components are working at level $j$ or bigger. Let $p_{1,1}=0.8, p_{1,2}=0.75, p_{2,1}=0.9, p_{2,2}=0.8$ and $p_{3,1}=0.75, p_{3,2}=0.7$.

  The only minimal $1$-working state of $\Sc$ is $(1,1,1)$ and the only minimal $2$-working state is $(2,2,2)$. The minimal $1$-failure states are $(2,2,0)$, $(2,0,2), (0,2,2)$ and the minimal $2$-failure states are $(2,2,1), (2,1,2)$, $(1,2,2)$. The $1$-reliability of the system is $R_1(\Sc)=0.54$ and $U_1(\Sc)=0.46$; the $2$-reliability is $R_2(\Sc)=0.42$ and $U_2(\Sc)=0.58$.
\end{Example}

In {\em cut systems} the focus is on the failure of the system, the exact counterpart of path systems. In this case the levels of the system, $\{0,\dots,M\}$, indicate growing levels of failure. Hence being at level~$0$ indicates that the system is completely functional and level $j>i$ indicates that the system is performing worse  at level~$j$ than at level~$i$ (\ie~the higher the level, the higher the intensity of failure). For each state $j$ of each of the components $c_i$ we denote by $q_{i,j}$ the probability that component $c_i$ is failing at  level $\ge j$. The $j$-unreliability $U_j(\Sc)$ of $\Sc$ is the probability that $\Sc$ is failing at level~$\ge j$; conversely, the $j$-reliability of $\Sc$ is $R_j(\Sc)=1-U_j(\Sc)$. A {\em cut system} is given at level $j$ by its set of $j$-failing states \ie~those tuples $(s_1,\dots,s_n) \in S_1\times\cdots\times S_n$ such that $\Phi(s_1,\dots,s_n)\geq j$. We say that a state $(s_1,\dots,s_n)\in S_1\times\cdots\times S_n$ is a {\em minimal $j$-failing state} or {\em minimal $j$-cut} if  $\Phi(s_1,\dots,s_n)\geq j$ and $\Phi(t_1,\dots,t_n)< j$ whenever all $t_i\leq s_i$ and at least in one case the inequality is strict. We say that a state $(s_1,\dots,s_n)\in S_1\times\cdots\times S_n$ is a {\em minimal $j$-working state} or {\em minimal $j$-path} if  $\Phi(s_1,\dots,s_n)< j$ and $\Phi(t_1,\dots,t_n)\geq j$ whenever all $t_i\geq s_i$ and at least one of the inequalities is strict. Cut systems are usually denoted by :F (for {\em fail}) in the literature.

\begin{Example}\label{ex:parallelF}
  Let $\Sc$ be a cut system with $3$ components $c_1,c_2,c_3$ such that $S=S_1=S_2=S_3=\{0,1, 2\}$. The structure function of $\Sc$ is given by $\Phi(s_1,s_2,s_3)=\max\{s_1,s_2,s_3\}$ \ie~$\Sc$ is a multi-state {\em parallel:F} system, the system fails at level $j>0$ or bigger whenever any of its components is failing at level $j$ or more. Let $q_{1,1}=0.25, q_{1,2}=0.2, q_{2,1}=0.2, q_{2,2}=0.1$ and $q_{3,1}=0.3, q_{3,2}=0.25$.

  The minimal $1$-failing states of $\Sc$ are $(1,0,0), (0,1,0), (0,0,1)$ and the minimal $2$-failing states of this system are $(2,0,0), (0,2,0), (0,0,2)$. The $1$-unreliability of the system is $U_1(\Sc)=0.58$ and $R_1(\Sc)=0.42$; the $2$-unreliability of the system is $U_2(\Sc)=0.46$ and $R_2(\Sc)=0.54$.	
\end{Example}

For clarity, in this paper we refer mainly to path systems. All the computations and considerations may be correspondingly applied to cut systems unless otherwise stated. We assume that the working or failure probabilities of the components in any system are independently distributed, although the method can also be applied to systems with dependent components. Also, we will consider that the probabilities of the systems' components to be in their different levels are constant in time (or equivalently, we consider reliability at a given instant $t$ or for the steady state of the system). The algebraic method can be applied to components' probabilities that vary in time, and to repairable systems and renewal processes. In those cases, we must consider the $j$-reliability of the system, \ie~the probability that the system is performing at level $j$ or better, and also the $j$-availability of the system at a given instant $t$ or an interval $\mathcal{T}$ of time, \ie~the probability that the system is performing at level $j$ or better for all $t\in \mathcal{T}$. Both notions have been extensively studied in the literature of multi-state systems, e.g. \cite{EPS78}. For complete introductions to multi-state system reliability and methods see \cite{KZ03,N11,TB17}.

\subsection{Algebraic reliability}\label{sec:algebraicReliability}
The use of commutative algebra (in particular of monomial ideals) in system reliability started in \cite{GNW02, GW04} in a close relation to improvements in inclusion-exclusion formulas and Bonferroni bounds \cite{D03}.
The approach we follow
was developed in a series of papers, e.g. \cite{SW09,SW15,MPSW19}. The main idea is to associate an algebraic object to a coherent system and obtain information about the structure and reliability
of the system by investigating the properties of
the algebraic object. In this section we make a brief self-contained description of the algebraic concepts involved and refer the interested reader to the cited series of papers for full details and proofs.

Let $\Sc$ be a path system with $n$ components and let $j\in\{0,\dots,M\}$ be one of the levels of the system. Let $F_j(\Sc)$ be the set of $j$-working states of the system and $\overline{F}_j(\Sc)$ the subset of minimal $j$-working states. Let us consider $P=\kb[x_1,\dots,x_n]$ a polynomial ring in $n$ indeterminates, one for each component of $\Sc$; here $\kb$ denotes any field of characteristic~$0$, but for clarity we assume that our coefficients are in $\mathbb{Q}$ or $\mathbb{R}$. To each state $s=(s_1,\dots,s_n)\in S_1\times\cdots\times S_n$ of $\Sc$ we associate the monomial $x^s=x_1^{s_1}\cdots x_n^{s_n}\in P$.

We denote by $pr(x^s)=\prod_{i=1}^n p_{i,s_i}$ the probability that the system is in a state $\ge s$.
In algebraic terms, having a state $t\ge s$ is equivalent to saying that the monomial $x^t$ is a multiple of $x^s$, \ie~$x^t$ is in $I=\langle x^s\rangle$, the ideal in $P$ generated by the monomial~$x^s$.
Now we consider the probability that the system is in a state greater then or equal to at least one of the states in $\{\mu_1,\dots,\mu_r\}$: this situation algebraically corresponds to the set of monomials $x^t$ belonging to the ideal  $I=\langle x^{\mu_1},\dots,x^{\mu_r}\rangle$.
  Thus we denote its probability by  $pr(I)=pr(\bigcup_{i=1}^r \langle x^{\mu_i}\rangle)$.

  The ideal generated by the $j$-working states of $\Sc$ is denoted by $I_j(\Sc)$ and is called the {\em $j$-reliability ideal} of $\Sc$. For a monomial ideal there is a unique minimal monomial generating set, denoted ${\rm MinGens}(I)$. Thus we observe that, due to the coherence property of $\Sc$, we have that ${\rm MinGens}(I_j(\Sc))$ is the set of the monomials corresponding to the minimal $j$-paths of $\Sc$
\[
 I_j(\Sc)=\langle x^{\mu}\mid \mu\in \overline{F}_j(\Sc) \rangle.
 \]
 
 From these definitions we have that the reliability of $\Sc$ is given by
 \[
R_j(\Sc)=pr(I_j(\Sc)).
\]

\begin{Remark}
The definition and description of the structure function of a multi-state coherent systems can be based on the sets of minimal $j$-paths or minimal $j$-cuts, which extend the notion of minimal paths and minimal cuts from binary systems \cite{BK94}. These special sets of state vectors are described as lower boundary points and upper boundary points to level $j$ in \cite{BK94}. Both sets can be formulated in algebraic terms as the minimal generators of the $j$-reliability ideal (lower boundary points) and maximal standard pairs (upper boundary points). For full details and a complete proof of this correspondence, see \cite{PSW20}.
\end{Remark}

Since $pr(I_j(\Sc))$ is expressed as the probability of a union, a natural choice for this computation is using the inclusion-exclusion principle, which in this case can be expressed as
\begin{equation}\label{eq:incl-excl}
 pr(I_j(\Sc))=\sum_{i=1}^r(-1)^{i+1}\sum_{\vert\sigma\vert=i}pr(\lcm(x^{\mu_s}\vert s\in\sigma)),
\end{equation}
 where $\sigma$ denotes subsets of $\{1,\dots,r\}$ and $\lcm$ denotes the least common multiple.

 A compact form of Equation (\ref{eq:incl-excl}) can be obtained by the multigraded Hilbert series of $I_j(\Sc)$. The multigraded Hilbert series of an ideal or module is a very important invariant in commutative algebra and algebraic geometry \cite{E95}. It is useful in our context because it provides a compact way to enumerate all monomials in a monomial ideal. The multigraded Hilbert series of an ideal $I\in P$, given by
 \[
 H_I(x_1,\dots,x_n)=\sum_{\mu\in\mathbb{N}^n}[x^{\mu}\in I]x^\mu,
   \]
   where the symbol $[x^\mu\in I]$ is equal to $1$ if $x^\mu$ is in $I$ and $0$ otherwise.
   The multigraded Hilbert series is an element of the formal power series ring $\mathbb{Z}[[x_1,\dots,x_n]]$. Observe that in this ring we have the identity $\frac{1}{1-x_i}=1+x_i+x_i^2+\cdots$ and hence one way to enumerate all the monomials in $P$ is to consider the summands of $\prod_{i=1}^n\frac{1}{1-x_i}$. Therefore the Hilbert series of $P$ is given by $H_P(x_1,\dots,x_n)=\prod_{i=1}^n\frac{1}{1-x_i}$. Just by multiplying every monomial by a given $x^\mu\in P$ one obtains that  $H_{\langle x^\mu\rangle}(x_1,\dots,x_n)=\prod_{i=1}^n\frac{x^\mu}{1-x_i}$. Now, since the set of monomials in a monomial ideal $I$ generated by $\{x^{\mu_1},\dots,x^{\mu_r}\}$ is the union of the sets of monomials in each of the ideals $\langle x^{\mu_i}\rangle$ then we have that
   \begin{equation}\label{eq:incl-excl-hilbert}
      H_I(x_1,\dots,x_n)=\sum_{i=1}^r(-1)^{i+1}\sum_{\vert\sigma\vert=i}\frac{\lcm(x^{\mu_s}\vert s\in\sigma)}{\prod_{j=1}^n(1-x_j)},
   \end{equation}
   which is the algebraic version of Equation (\ref{eq:incl-excl}). Let $HN_I(x_1,\dots,x_n)$ denote the numerator of the Hilbert series of the ideal $I$ and let $\Sc$ be a coherent path system as in Section \ref{sec:coherentSystems}. Let $pr(HN_{I_j(\Sc)}(x_1,\dots,x_n))$ denote the formal substitution of every $x^\mu$ by $pr(x^\mu)$ in the numerator of the multigraded Hilbert series of $I_j(\Sc)$. The direct relation between equations (\ref{eq:incl-excl}) and (\ref{eq:incl-excl-hilbert}) allows us to establish the fundamental identity of the algebraic approach to system reliability
   \begin{equation}\label{eq:reliability-hilbert}
     R_j(\Sc)=pr(I_j(\Sc))=pr(HN_{I_j(\Sc)}(x_1,\dots,x_n)).
   \end{equation}
   Hence any way to obtain $HN_{I_j(\Sc)}(x_1,\dots,x_n)$ gives us a way to compute $R_j(\Sc)$. Of course a direct one, although very redundant in general, is Equation (\ref{eq:incl-excl-hilbert}) by means of the inclusion-exclusion principle. Other more efficient and compact ways to obtain $HN_{I_j(\Sc)}(x_1,\dots,x_n)$ are described in \cite{B97}.

   An important feature of the inclusion-exclusion formulas is that they can be truncated to obtain the so called Bonferroni bounds \cite{D03}. More precisely, we have that
   \begin{equation}
     \begin{aligned}\label{eq:incl-excl-bounds}
       pr(I_j(\Sc))\leq\sum_{i=1}^t(-1)^{i+1}\sum_{\vert\sigma\vert=i}pr(\lcm(x^{\mu_s}\vert s\in\sigma))\mbox{ for } t\leq r \mbox{ odd,}\\
       pr(I_j(\Sc))\geq\sum_{i=1}^t(-1)^{i+1}\sum_{\vert\sigma\vert=i}pr(\lcm(x^{\mu_s}\vert s\in\sigma)) \mbox{ for } t\leq r \mbox{ even.}
     \end{aligned}   
   \end{equation}
   One way to obtain the multigraded Hilbert series of a monomial ideal~$I$ is by constructing a multigraded free resolution of $I$ and read $HN(I)$ from the data in the resolution. Every ideal $I\subseteq P=\kb[x_1,\dots,x_n]$ can be described as a module in terms of what is called a {\em free resolution}, which is a series of free modules and morphisms among them. A free module is a direct sum of copies of $P$ with the usual grading shifted by some degree $d\in\mathbb{N}$ denoted by $P(-d)$. In the case of monomial ideals we can also have multigraded resolutions, in which the degree shifts are given by multidegrees $(d_1,\dots,d_n)\in\mathbb{N}^n$ and the shifted copies of $P$ are denoted by $P(\-\mu)$. A multigraded free resolution of a monomial ideal $I$ is of the form
   \[
0\longrightarrow \bigoplus_{j=1}^{r_d}P(-\mu_{d,j})\stackrel{\partial_d}{\longrightarrow} \cdots \stackrel{\partial_2}{\longrightarrow}\bigoplus_{j=1}^{r_1}P(-\mu_{1,j})\stackrel{\partial_1}{\longrightarrow}P/I\longrightarrow 0,
\]
where the $\partial_i$ are graded module morphisms (here called {\em differentials}), $d$ is the length of the resolution, the $r_i$ are called {\em ranks} of the resolution and the $\mu_{i,j}$ for each $i$ are the multidegrees of the $i$-th module of the resolution. Given an ideal $I$ one can build different resolutions and among them there is a distinguished one called the {\em minimal free resolution} which is unique up to isomorphisms and is characterized by having smallest ranks among all the possible resolutions of $I$. The ranks of the minimal free resolution of $I$ are called the {\em Betti numbers} of $I$ and are a fundamental invariant of $I$ \cite{E95}.

Now, given any multigraded free resolution of a monomial ideal $I$ we have the following expression for $HN_I(x_1,\dots,x_n)$
\begin{equation}
HN_I(x_1,\dots,x_n)=\sum_{i=1}^d(-1)^i\sum_{j=1}^{r_i}x^{\mu_{i,j}}.
\end{equation}
This expression can be truncated as in (\ref{eq:incl-excl-bounds}) and produces the following bounds for $R_j(\Sc)$, see \cite{SW09}
 \begin{equation}
     \begin{aligned}\label{eq:resolutionBounds}
       R_j(\Sc)\leq\sum_{i=1}^t(-1)^{i+1}\sum_{j=1}^{r_i}pr(x^{\mu_{i,j}})\mbox{ for } t\leq r \mbox{ odd,}\\
       R_j(\Sc)\geq\sum_{i=1}^t(-1)^{i+1}\sum_{j=1}^{r_i}pr(x^{\mu_{i,j}})\mbox{ for } t\leq r \mbox{ even.}\\
     \end{aligned}   
 \end{equation}
 Among this type of bounds, those given by the minimal multigraded free resolution of $I_j(\Sc)$ are the tightest, cf. \cite{SW09}.
 
 \begin{Example}
   Consider the series:G system $\Sc$ studied in Example ~\ref{ex:seriesG}. The $j$-reliability ideals of $\Sc$ are $I_1(\Sc)=\langle x_1x_2x_3\rangle$ and $I_2(\Sc)=\langle x_1^2x_2^2x_3^2\rangle$. The minimal free resolution of $I_1(\Sc)$ has length $1$ and the only free module of this resolution has multidegree $(1,1,1)$, hence $R_1(\Sc)=pr(x_1x_2x_3)=0.54$. An equivalent computation shows that $R_2(\Sc)=pr(x_1^2x_2^2x_3^2)=0.42$.
 \end{Example}

 \begin{Example}
   Consider now the parallel:F system $\Sc$ in Example \ref{ex:parallelF}. The $j$-unreliability ideals of $\Sc$ are $I_1(\Sc)=\langle x_1,x_2,x_3\rangle$ and $I_2(\Sc)=\langle x_1^2,x_2^2,x_3^2\rangle$. The minimal free resolution of $I_2(\Sc)$ has length $3$ and the multidegrees of its modules are $\mu_{1,1}=(2,0,0), \mu_{1,2}=(0,2,0), \mu_{1,3}=(0,0,2), \mu_{2,1}=(2,2,0), \mu_{2,2}=(2,0,2), \mu_{2,3}=(0,2,2)$ and $\mu_{3,1}=(2,2,2)$ hence the $2$-unreliability of $\Sc$ is given by
   \begin{align*}
     U_2(\Sc)&=pr(x_1^2)+pr(x_2^2)+pr(x_3^2)-(pr(x_1^2x_2^2)+pr(x_1^2x_3^2)+pr(x_2^2x_3^2))\\
     &+pr(x_1^2x_2^2x_3^2)=0.55-0.095+0.005=0.46.
   \end{align*}
   Observe that truncating this expression we obtain a first upper bound of $0.55$ and a first lower bound of $0.455$ for $U_2(\Sc)$.
 \end{Example}
 
The algebraic methodology for reliability computation using monomial ideals is based on two main principles. The first one is to avoid as much redundancy as possible when enumerating the states needed for the final reliability computation. This is provided by the possibility of using different resolutions to express the numerator of the Hilbert series of the system's ideals. In this respect, a fast computation of the minimal resolution or close-to-minimal resolutions is the main component of our approach. The second principle is that this methodology can be approached as a recursive procedure, computing the Hilbert series of an ideal in terms of the Hilbert series of smaller ideals. Recursion is usually very efficient in reliability computations and is used in other methodologies, such as the UGF method \cite{L05}, the decomposition and factoring methods \cite{KZ03,TB17} or ad-hoc methods for particular systems, see \cite{MLAD15} for instance.

To achieve the aforementioned two principles, we use in our algorithm the Hilbert series expression given by ranks of resolutions computed by Mayer-Vietoris trees, described in \cite{S09}. They are a fundamental tool in our implementation of the algebraic method for system reliability, for they are a fast and efficient recursive algorithm and their output is a resolution that is in most cases minimal or very close to minimal. This is the main ingredient of the efficiency of the class described in Section \ref{sec:class} which is demonstrated in some examples in Section \ref{sec:examples}.

\begin{Remark}
  The principles of avoiding redundant computations and using recurrence are the keys of the UGF method, which also uses an algebraic formalism. They are however used in a slightly different way compared to our approach. In the first place, the UGF method avoids the redundancy by collecting like terms during the recursive computation. In our approach we reduce redundant terms by checking divisibility of the involved monomials during the computation and hence not using terms that would be cancelled in the final evaluation of the expression. With respect to the recursive methodology, the recursions used in the UGF method needs a certain structure for the base cases (e.g. series or parallel systems as base case to stop the recursion) and the recursion is applied by blocks or subsystems. In our case the recursion is applied by selecting a single generator of the ideal in each step. This reduces the need of previous knowledge of the system's structure and is therefore very general, but has the disadvantage of eventually producing larger recursion trees. 

These observations suggest that both methods could mutually benefit from considering the stronger points of the other to improve their own performance. A full comparison is certainly worth considering, but it is beyond the scope of this paper.
\end{Remark}
 
 \begin{Example}\label{ex:doubleBridge}
Consider the double bridge binary system in Figure \ref{network}, taken from \cite{D03}. It has $5$ nodes and $8$ connections and we assume that only its connections are subject to failure. The minimal cuts of the system are
$123$, $1258$, $13457$, $1478$, $2346$, $24568$, $3567$ and $678$. The cut ideal
corresponding to this network is then
{\small 
\begin{align*}
I=&\langle
  x_1x_2x_3,\; x_1x_2x_5x_8,\; x_1x_3x_4x_5x_7,\; x_1x_4x_7x_8,\; x_2x_3x_4x_6,\\
 &x_2x_4x_5x_6x_8,\; x_3x_5x_6x_7,\; x_6x_7x_8\rangle.
\end{align*}
}

\begin{figure}
\begin{center}
\includegraphics{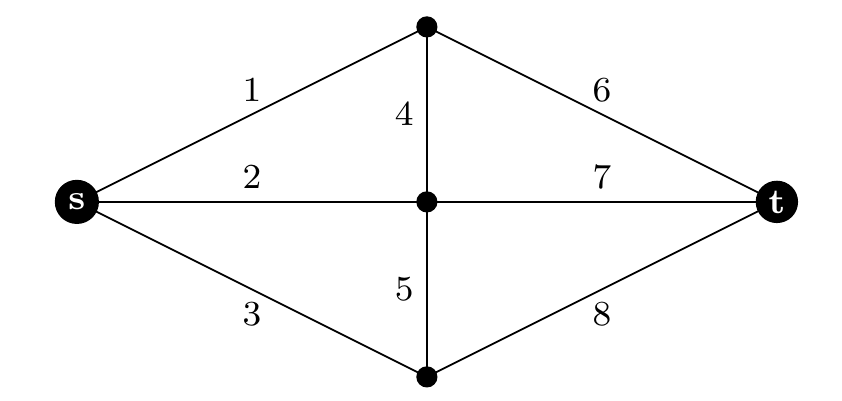}
\caption{Double bridge network.}\label{network}
\end{center}
\end{figure}

To compute the (un)reliability of this system we can use the Hilbert series of its cut ideal, which has $8$ generators. Using the Taylor resolution corresponds to the inclusion-exclusion method, and uses in this case $255$ summands for the Hilbert series. Another resolution is the Scarf resolution, that is equivalent to apply the method of \emph{abstract tubes} \cite{D03,GW04}, which in this case is much less redundant, giving an expression for the Hilbert series that uses only $51$ summands. Using the minimal free resolution gives an expression using $43$ summands and is the most compact form achievable with this methodology. If we use the expression given by the inclusion-exclusion method and Bonferroni bounds obtained by truncation, we have that the eighth bound is sharp, while in the case of the Scarf or minimal resolutions already the fourth bound is sharp.
\end{Example}

\begin{Remark}
The network in Example \ref{ex:doubleBridge} can be handled efficiently using the factoring algorithm. The domination invariant \cite{SC83} of this system  is $4$, which implies that the reliability function can be expressed very easily as the sum of just four terms. Observe that this is the same number of terms in the minimal resolution of the reliability ideal of this system.
\end{Remark}
 
 \subsection{Duality}
 Given a structure function $\Phi$ its dual $\Phi^D$ with respect to $\tb\in\NN^n$ is given by (cf. \cite{EPS78})

 \begin{equation}\label{eq:systemDualityBin}
\Phi^D(s_1,\dots,s_n)=M-\Phi((t_1-s_1,\dots,t_n-s_n)).
\end{equation}

\begin{Example}\label{ex:dual1}
Consider a binary series:G system $\Sc$ with three components where $\Phi(s_1,s_2,s_3)=\min\{s_1,s_2,s_3\}$. We have $\Phi^D(s_1,s_2,s_3)=0$ if and only if $(s_1,s_2,s_3)=(0,0,0)$ hence the minimal working states of the dual system are $(1,0,0),(0,1,0)$ and $(0,0,1)$, which correspond to a parallel system. The dual of a series system is always a parallel system and vice-versa.
\end{Example}

There is a notion of duality in monomial ideals, called {\em Alexander duality} \cite{MS04}. To describe it we use the following notation. Given a vector $\mu\in\NN^n$, we denote by $\mathfrak{m}^\mu$ the monomial ideal
\[
\mathfrak{m}^\mu=\langle x_i^{\mu_i}\mid \mu_i\geq 1\rangle.
\]
Given two vectors $\mu$ and $\nu$ in $\NN^n$ let $\mu\setminus\nu$ the vector whose $i$'th coordinate is $\mu_i+1-\nu_i$ if $\nu_i\geq 1$ and $0$ otherwise. 
\begin{Definition}
  Let $I\subset \kb[x_1,\dots,x_n]$ be a monomial ideal, ${\rm MinGens}(I)$ its minimal set of monomial generators, and $x^\nu=\lcm({\rm MinGens}(I))$. The Alexander dual of I is the intersection
    \[
I^D=\bigcap_{x^\mu\in{\rm MinGens}(I)}\mathfrak{m}^{\nu\setminus\mu},
    \]
where $\mathfrak{m}^{(s_1,\dots,s_n)}$ denotes the monomial ideal $\langle x_i^{s_i}\mid s_i\geq 1\rangle$ 
  \end{Definition}

\begin{Remark}
  Given a coherent system $\Sc$ its $j$-reliability ideal $I_j(\Sc)$ is generated by the monomials corresponding to  its minimal $j$-paths.
  The ideal of its dual system $I_j(\Sc^D)$
  is generated by the monomials corresponding to maximal $j$-cuts of $\Sc$ and may be seen as the ideal generated by the maximal standard pairs of $I_j(\Sc)$ \cite{PSW20}. These can be computed using the Alexander dual of the artinian ideal $I_j(\Sc)+\langle x_i^{M_1+1},\dots,x_n^{M_n+1}\rangle$ \cite{BS09}.
\end{Remark}

We can use the dual ideal of a system to compute its reliability in the following way. Let $\overline{pr}(x^\mu)=\prod_{i=1}^n(1-p_{i,\mu_i+1})$ \ie~the product of the probabilities that each component $i$ is in a state less than or equal to $\mu_i$. We denote by $\nu=(M_1,\dots,M_n)$ the vector of maximal possible levels of the components. Let $\overline{I}_{j}(\Sc)$ be the ideal generated by the monomials $\{x^{\overline{\mu}}\mid x^\mu \mbox{ is a generator of }I_j(\Sc)\}$. We consider the ideal $\overline{I}_j(\Sc)^D$ and compute $H_{\overline{I}_j(\Sc)^D}(x_1,\dots,x_n)$. We obtain $U_j(\Sc)=1-R_j(\Sc)$ by formally substituting each monomial $x^{\mu}$ in $HN_{\overline{I}_j(\Sc)^D}(x_1,\dots,x_n)$ by $\overline{pr}(\frac{x^\nu}{x^\mu})$.

%

\begin{Example}
  Consider the system in Example \ref{ex:dual1}. We have that $I_1(\Sc)=\langle x_1x_2x_3\rangle$, then $I_1(\Sc)^D=\langle x_1,x_2,x_3\rangle$ and $\overline{I}_1(\Sc)^D=\langle x_1,x_2,x_3\rangle$. By using the minimal free resolution of $\overline{I}_1(\Sc)^D=\langle x_1,x_2,x_3\rangle$ we have that $HN_{\overline{I}_1(\Sc)^D}(x_1,\dots,x_n)=(x_1+x_2+x_3)-(x_2x_3+x_1x_3+x_1x_2)+x_1x_2x_3$. Hence, if we set the probabilities $p_{1,1}=0.8,\, p_{2,1}=0.9$ and $p_{3,1}=0.75$, we get
  \begin{align*}
    U_1(\Sc)&=\overline{pr}(x_1x_2)+\overline{pr}(x_1x_3)+\overline{pr}(x_2x_3)-\overline{pr}(x_1)-\overline{pr}(x_2)-\overline{pr}(x_3)+\overline{pr}(1)\\
          &=0.25+0.1+0.2-(0.025+0.05+0.02)+0.005=0.46,\\
  \end{align*}
 and we obtain $R_1(\Sc)=0.54$. Observe that in the equality above, $\overline{pr}(1)=\overline{pr}(x_1^0x_2^0x_3^0)=pr(x_1\leq 0)pr(x_2\leq 0)pr(x_3\leq 0)=0.005$.
\end{Example}

\section{Algebraic reliability class in \texttt{CoCoALib}}\label{sec:class}
The good performance of an algorithm depends also on the efficiency of its implementation. In this section we give the interested reader some technical details on the implementation of our algorithms and some of the decisions we made, like the choice of data types and the structure of the algorithms. These decisions contribute to the actual performance of the algorithms in terms of memory usage and CPU time. Also, we describe the \verb|CoCoALib| library, which provides convenient implementations of the main algebraic structures we need. We hope these descriptions, although not fully detailed, make it easier for engineers and reliability practitioners to practically use these algorithms or incorporate them into their own software, and also make it easy to reproduce our results, experiments and benchmarks.

\subsection{CoCoALib}
\verb|CoCoALib|, for Computations in Commutative Algebra, is an open source \verb|C++| software library principally based on multi- variate monomials and polynomials and devoted to algebraic geometry. It is the computational core of the CoCoA software system \cite{CoCoA}.
A crucial aspect of CoCoALib is that it was designed from
the outset to be an open-source software library.  This initial decision, together with the desire to help the
software prosper, has many implications: {\em e.g.} designing a
particularly clean interface for all functions with comprehensive
documentation.  This cleanliness makes it easy to integrate CoCoALib
into other software in a trouble-free manner. The library is fully documented, and also comes with about 100 illustrative example programs.
CoCoALib reports errors using \verb|C++| exceptions, while the
library itself is exception-safe and thread-safe.  The
current source code follows the \verb|C++14| standard. The main features of the design of CoCoALib are:
\begin{itemize}
\item it is well-documented, free and open source \verb|C++| code (under the GPL~v.3 licence);

\item the design is inspired by, and respects, the underlying mathematical structures;

\item the source code is clean and portable;

\item the user function interface is natural for mathematicians, and easy to memorize;

\item execution speed is good with robust error detection.
\end{itemize}

The design of the library (and its openness) was chosen
to facilitate and encourage ``outsiders'' to contribute.  There are
two categories of contribution: code written specifically to become part
of CoCoALib, and stand-alone code written without considering its
integration into CoCoALib. The library has combined some of the features of various external libraries
into CoCoALib. Such as \textit{Frobby}
(see~\cite{Frobby}) which is specialized for operations on monomial
ideals. Other integrations are with \textit{Normaliz} library for computing with affine monoids or
rational cones and \textit{GFanLib} which is a \verb|C++| software library for computing Gr\"obner fans and tropical varieties.

\subsection{The class description}
We have implemented within CoCoALib a set of \verb|C++| classes for making computations in algebraic reliability. The UML class diagram is depicted in Figure \ref{fig:classDiagram} in the Appendix. Our main class is the {\em abstract class}
\verb|CoherentSystem| which consists of a series of levels and a matrix of probabilities. The levels are stored in a \verb|std::vector| (an efficient structure of the \verb|C++| language) in which each component is an instance of the class \verb|CoherentSystemLevel|, and the probabilities are given by a \verb|vector| of \verb|vectors| of type \verb|double| where the $j$'th entry of the $i$'th vector corresponds to $p_{i,j}=p(s_i\geq j)$, the probability that the level of the $i$'th component of the system is bigger than or equal to $j$. Each instance of the class \verb|CoherentSystemLevel| consists basically of an ideal and its dual, which are objects of the CoCoALib class \verb|ideal|. Also, we store as member fields their Mayer-Vietoris trees, which play the role of multigraded free resolutions optimized for monomial ideals.

The {\em concrete classes} inheriting from the class \verb|CoherentSystem| are \verb|CoherentSystemPath| and \verb|CoherentSystemCuts| which respectively represent :G systems in which the levels and probabilities denote working states, and :F systems in which the levels and probabilities represent failures, as seen in Section ~\ref{sec:coherentSystems}.
For any instance of these two concrete classes, and hence of the abstract class \verb|CoherentSystem| we can call the following member functions:
\begin{description}
\item[myMinimalPaths] Receives a level and gives a vector of vectors of type \verb|long|. Each of these vectors is a minimal path of the system at the given level.
\item[myMinimalCuts] Receives a level and gives a vector of vectors of type \verb|long|. Each of these vectors is a minimal cut of the system at the given level.
\item[myReliability] Receives a level $j$ and computes $R_j(\Sc)$.
\item[myUnreliability] Receives a level $j$ and computes $U_j(\Sc)$.
\item[myReliabilityBounds] Receives a level $j$ and computes bounds for $R_j(\Sc)$ given by the resolution obtained by the Mayer-Vietoris tree of $I_j(\Sc)$ as given in Equation (\ref{eq:resolutionBounds}).
\item[myUnreliabilityBounds] Receives a level $j$ and computes bounds for $U_j(\Sc)$ given by the resolution obtained by the Mayer-Vietoris tree of $I_j(\Sc)$ computed from the bounds for $R_j(\Sc)$ given in Equation (\ref{eq:resolutionBounds}).  
\end{description}

In addition, for :G systems given by its sets of paths, we have implemented two more bounds, described by G\r{a}semyr and Natvig in \cite{GN17}:
\begin{description}
\item[GNMaxMinPathBound]{Let $\bf{y}^m$, $m=1,\dots,M_p$ the minimal paths of $\Sc$ for level $j$, the following lower bound for $R_j(\Sc)$ is given in \cite{FN85}:
  \[
    l'^{j}({\bf p})=\max_{1\leq m\leq M_p}\left( \prod_{i=1}^n p_i^{y_i^m}\right) 
  \]
}
\item[GNCoproductMinCutsBound]{Let $\bf{z}^m$, $m=1,\dots,M_c$ the set of minimal cut vectors of $\Sc$ for level $j$, then we have the following lower minimal bound for $R_j(\Sc)$ \cite{FN85}:
    \[
       l^{**j}({\bf p})=\prod_{m=1}^{M_c}\coprod_{i=1}^np_i^{z_i^m+1}
    \]
    where for $p_i\in[0,1]$ we define $\coprod_{i=1}^n=1-\prod_{i=1}^n(1-p_i)$.
    }
  \end{description}

  When computing the functions \verb|myReliability|, \verb|myUnreliability|, \verb|myReliabilityBounds| or \verb|myUnreliabilityBounds| the object checks its ideal and its dual ideal, and chooses whichever of them has a smallest number of minimal generators to perform the actual computation. To compute duals of ideals we use the \verb|Frobby| library, in particular the function \verb|FrbAlexanderDual| which is in general a fast computation. Once the ideal is chosen, we check whether the system has already computed its Mayer-Vietoris tree. If it is not yet computed, it is computed and stored in the corresponding class member field. Then the Mayer-Vietoris tree is used to retrieve the required value or bounds for reliability or unreliability.

  \section{Examples of use}\label{sec:examples}
  In this section we apply our \verb|C++| class to some examples of reliability computations. We use binary networks and multi-state systems. We consider systems in which their components have independent identically distributed probabilities as well as systems in which the components' probabilities are independent but not identically distributed. All the computations in this section have been implemented by the authors and executed in an HP Z-book laptop\footnote{CPU: intel i7-4810MQ, 2.80 GHz. RAM: 16Gb}.
  \subsection{Test examples}
  First, we validate our algorithms with a set of diverse examples of multi-state systems found in the literature. We selected systems of different nature so that we can test our algorithms with examples featuring different characteristics. Table \ref{table:validation} shows the results of these tests. The first column of the table indicates the name of each example (see description below), $n$ indicates the number of variables and $M$ the number of levels of the system (not counting the complete failure level or level $0$). Column $M_i$ indicates the number of levels of each component and column $gens(I_j)$ indicates the number of minimal generators of the $j$-reliability ideal for each level $j=1,\dots,M$. The set of test examples consists of the following:
  \begin{itemize}
  \item[-]{\texttt{Army Battle Plan} is taken from the classical paper \cite{BK94}. It is a customer-driven multi-state system with $5$ different states and $4$ components (two binary components and two three-level components), the probabilities of the different components are independent but not identical.}
  \item[-]{\texttt{Bin.S-P} is a binary Series-Parallel system taken from \cite{LL03} (Example 4.5) which has seven independent not identical components and two levels. }
  \item[-] {\texttt{MAX+MIN,TIMES} is a multi-state system with $5$ components and $7$ levels whose structure function is given by 
   $$\Phi(x_1,\dots,x_5)=\left( \max\{x_1,x_2\}+\min\{ x_3,x_4\}\right)\times x_5,$$
   and the details on components' and system's levels and probabilities (not i.i.d) are given in \cite{LFD10}, Example 4.7.}
  \item[-]{\texttt{Bridge Flow Network} is a multi-state network with $5$ edges with different weights considered as flows. The states of the system are given by the possible flows through the network. The example considers the probability of a total flow of at least three units (\ie~the system is in level $j=3$). The details on states and probabilities are given in \cite{SXD10} Example 4, see also \cite{TB17} Example 5.14.}
 \item[-]{\texttt{Dominant MS binary-imaged} system is a multi-state system with three i.i.d. components. Both the system and components have four different states. It is presented as Example 12.21 in \cite{KZ03} to illustrate the concept of multi-state dominant binary-imaged system.}
 \item[-]{\texttt{MS Cons.k-out-of-n} is a multi-state consecutive $k$-out-of-$n$ system with $3$ components and three levels. It is example 12.18 in \cite{KZ03}.}
  \end{itemize}
  
  \begin{table}
  \begin{tabular}{llcll}
    \hline
    {Example}&{$n$}&{$M$}&{$M_i$}&$gens(I_j)$\\
      \hline
       \hline
    \texttt{Army Battle Plan}&4&4&2,2,3,3&2,4,5,5\\
    \texttt{Bin.S-P}&7&1&$1\;\forall i$&3\\
    \texttt{MAX+MIN,TIMES}&5&6&3,2,2,3,2&4,3,4,3,2,1\\
    \texttt{Bridge Flow Network}&5&3&3,1,2,1,2&3 for $j=3$\\
    \texttt{Dominant MS binary-imaged system}&3&3&3,3,3&3,2,1\\
    \texttt{MS Cons. k-out-of-n}&3&3&3,3,3&1,2,1\\


    \hline
    \end{tabular}\caption{Test examples of multi-state systems.}\label{table:validation}
  \end{table}
  
  \subsection{Source to terminal networks}
  One of the main problems in reliability engineering is Network Reliability, see for instance \cite{B79,TB17} for a comprehensive account and  \cite{GM20} for a recent algorithm. In this problem we consider a network in which one vertex is selected as {\em source vertex} and one or more vertices are selected as {\em target vertices}. Each of the connections in the network has a certain probability to be working, and the problem is to compute the probability that there exists at least one source-to-target path composed by operational connections. Usually the networks are binary \ie~the system and all of its components have only two possible levels, although the multi-state version has also been considered \cite{BTZ18,ZLB20}.
  
  \subsubsection{GARR: Italian Research and Education Network}
\begin{figure}
      \includegraphics[scale=0.4]{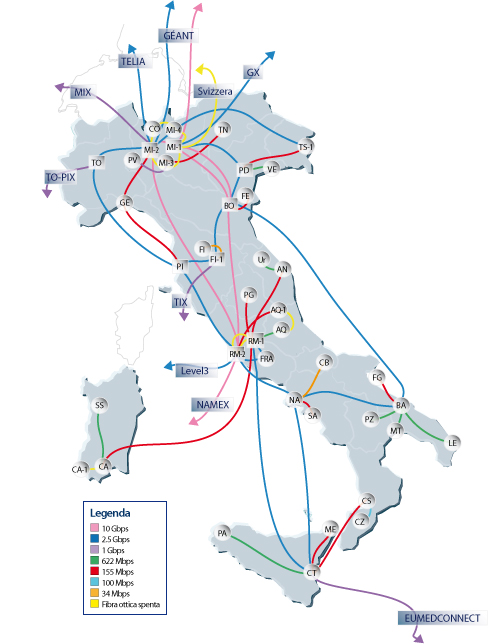}
  \caption{Map of the GARR network in 2008}\label{fig:garr}
\end{figure}

Our first example is the GARR Italian network. The motivation to use this example is to show the performance of our algorithms in a real-life system that has already been studied in the literature. Figure \ref{fig:garr} shows the official 2008 map of the backbone of the GARR network in Italy, which interconnects universities, research centers, libraries, museums, schools and other education, science, culture and innovation facilities, see \texttt{http//:www.garr.it}. The network was at the time formed by $41$ nodes and $52$ connections. Table \ref{table:garr} shows the results of some reliability computations in this network. First, we use \verb|TO| as source node and \verb|CT| as terminal node, and then we use \verb|TS1| as source node and \verb|NA| as terminal node. In both cases we consider an identical independent probability $p$ for all the connections and compute the source to terminal reliability for $p=0.9, 0.95$ and $0.99$. The last two rows in  Table \ref{table:garr} correspond to the same computations in \cite{TB17} (Example 5.7), observe that the differences are due to the fact that the authors in \cite{TB17} use a slightly different network which has $42$ nodes and hence some different connections and a different number of minimal paths in each example. Since our algorithms can also treat the case of non-identical probabilities, we assigned probability $0.99$ to all $10$Gbps connections ($4$ connections), $0.95$ to all $2.5$ Gbps connections ($14$ connections) and $0.9$ to the rest of the connections. The results are shown in the last column of the table. All our computations in this table took less than one second.

\begin{table}
  \begin{tabular}{llccccc}
    \hline
    \multirow[t]{2}{*}{S-node}&\multirow[t]{2}{*}{T-node}&\multirow[t]{2}{*}{\# Minpaths}&\multicolumn{3}{c}{i.i.d probabilities}& non i.i.d.\\

    \cline{4-6}
    \multicolumn{3}{l}{ }&0.9&0.95&0.99& \\
    \hline
    TO&CT&212&0.977344&0.994704&0.999798&0.994352\\
    TS1&NA&223&0.985203&0.996890&0.999895&0.993917\\
    \hdashline
    TO*&CT*&196&0.977428&0.994713&0.999797& --\\
    TS1*&NA*&168&0.975771&0.994486&0.999795& --\\  
    \hline
    \end{tabular}\caption{Reliability computations for the GARR 2008 network.}\label{table:garr}
  \end{table}

  \subsubsection{Random networks}
\begin{figure}
      \includegraphics[scale=0.45]{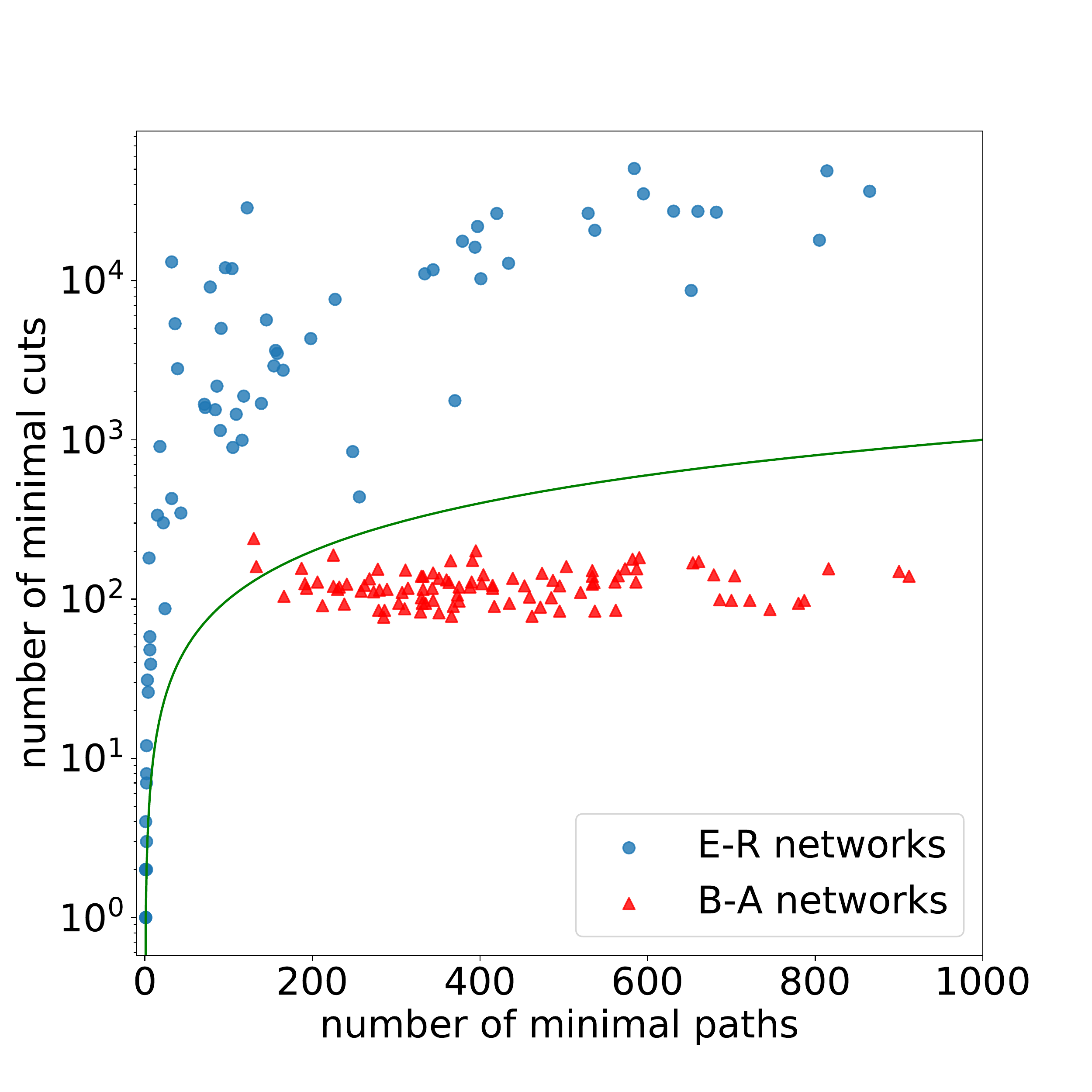}
  \caption{Number of minimal paths and minimal cuts in Erd\H{o}s-R\'enyi and Barabasi-Albert networks. The line indicates number of minimal paths equals number of minimal cuts.}\label{fig:paths_vs_cuts}
\end{figure}

\begin{figure}
     \resizebox{\textwidth}{!}{
\begin{tabular}{cc}
 \includegraphics[scale=0.8]{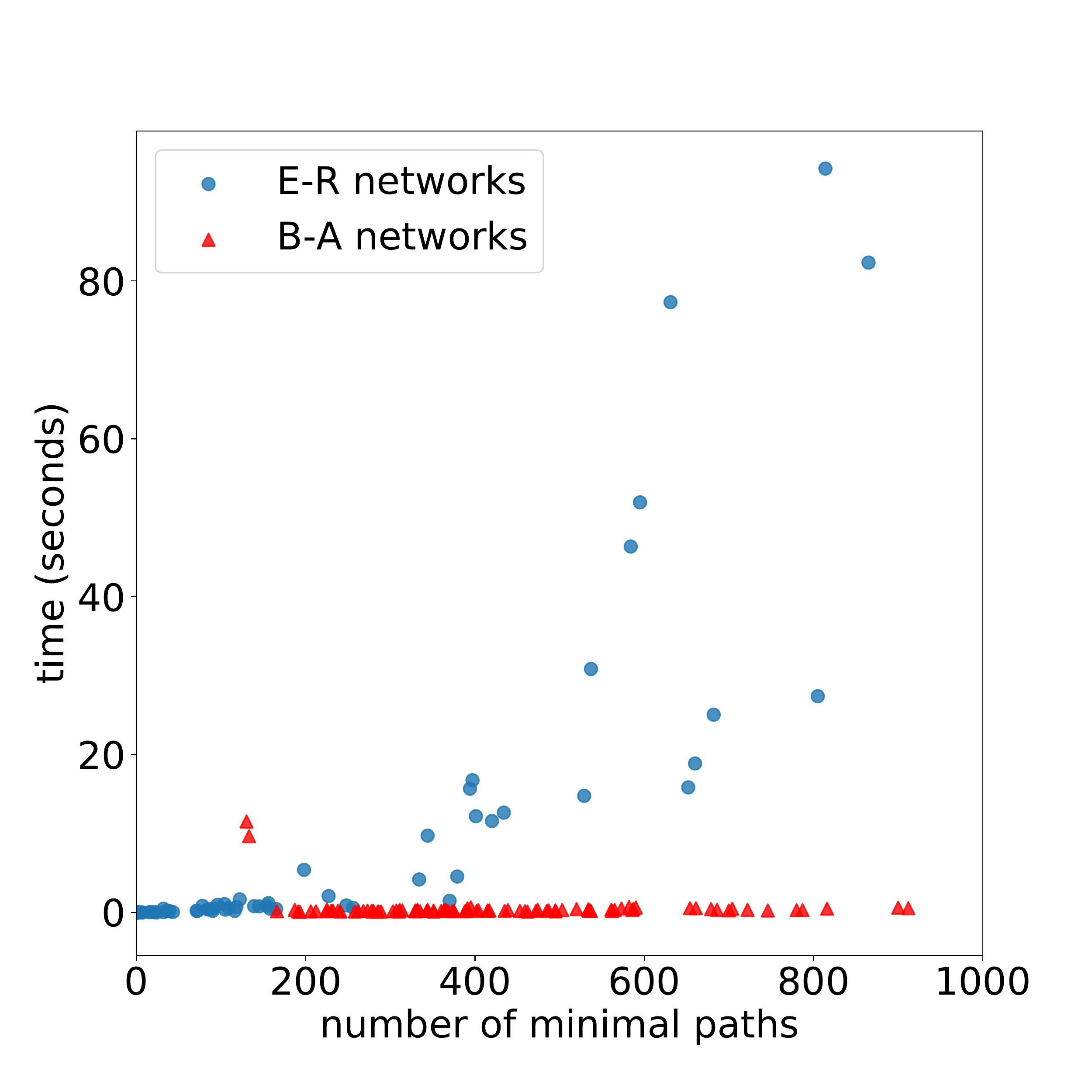} &   \includegraphics[scale=0.8]{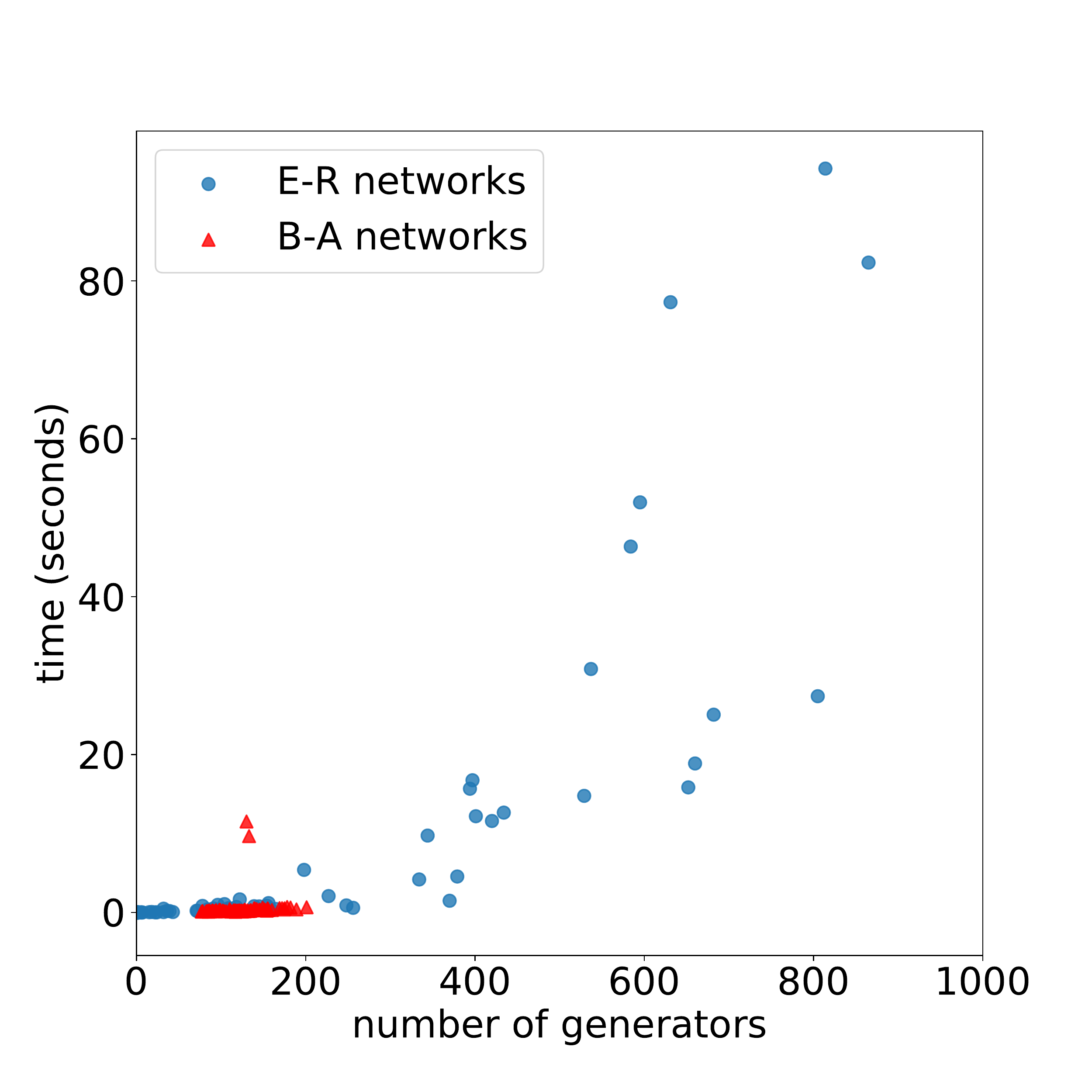}
 \\
(a) &(b) \\[6pt]
\end{tabular}}
\caption{Times for reliability computation on Erd\H{o}s-R\'enyi and Barabasi-Albert graphs.}
\label{fig:ERBA}
\end{figure}

 Our second example is a randomly generated set of networks. This is a convenient way to generate a big number of examples not having a regular structure (like for instance series-parallel systems or $k$-out-of-$n$ systems and variants), and therefore represents a good set of benchmarks for the application to general systems. We demonstrate our algorithms' performance in several random networks generated following the Erd\H{o}s-R\'enyi model $ER(n,p)$ \cite{ER59} and Barabasi-Albert model $BA(n,m)$ \cite{BA99}. These models generate networks with different characteristics such as degree distribution, modularity, etc. We compute the reliability of $100$ random Erd\H{o}s-R\'enyi networks with $n=40$, $p=0.05$ and $100$ Barabasi-Albert networks with $n=10$ and $m=4$; we chose randomly one source and one terminal node in each case. The number of minimal paths varies between $100$ and $1000$ in both cases. However, the relation between the number of minimal paths and minimal cuts is significantly different in the two types of networks. Erd\H{o}s-R\'enyi networks tend to have many more minimal cuts with respect to the number of minimal paths, while the situation is the opposite for Barabasi-Albert networks, see Figure \ref{fig:paths_vs_cuts}. In the case of Barabasi-Albert networks our algorithms compute the reliability of the network using the dual ideal, since it is smaller in most cases. The reliability of the Erd\H{o}s-R\'enyi examples was always computed using the minimal path ideal. Times for the computation of the reliability of these networks are shown in Figure~\ref{fig:ERBA}. The figures show that the times depend greatly on the number of minimal paths, but the topology of the network influences the algebraic characteristics of the ideals. Observe that there are two cases of Barabasi-Albert graphs in which the number of minimal paths is smaller than the number of minimal cuts and hence the path ideal was used for the computation, which results in higher computation times compared with the cases in which the dual was used. The resolutions of these networks ideals are much shorter in the dual case and hence the results. 

\subsection{Multi-state generalized $k$-out-of-$n$ systems}\label{subsec:kn}
Our final example is multi-state generalized $k$-out-of-$n$ systems. We include this example since they are one of the most important types of systems studied in the reliability engineering literature, both in their binary and multi-state versions. A binary $k$-out-of-$n$:G system is a system with $n$ components that is in a working state whenever at least $k$ of its components are working. The multi-state version of this kind of systems has received several definitions in the literature, see \cite{PSW20} for a review. A general definition is that of {\em generalized multi-state $k$-out-of-$n$ systems}, see \cite{HZW00}:

\begin{Definition}
  An $n$-component systems is called a generalized $k$-out-of-$n$:G system if $\phi(s_1,\dots,s_n)>j$, $1\leq j\leq M$ whenever there exists an integer value $l$, $(j\leq l\leq M)$ such that at least $k_l$ components are in state $l$ or above.
\end{Definition}

If we denote by $N_j$ the number of components of the system that are in state $j$ or above then this definition can be rephrased by saying that~$\phi(S)\geq j$ if
\begin{align*}
  N_j&\geq k_{j}\\
  N_{j+1}&\geq k_{j+1}\\
  \vdots&\\
  N_M&\geq k_M.\\
 \end{align*}

\begin{figure}[b]
     \resizebox{\textwidth}{!}{
\begin{tabular}{cc}
 \includegraphics[scale=0.5]{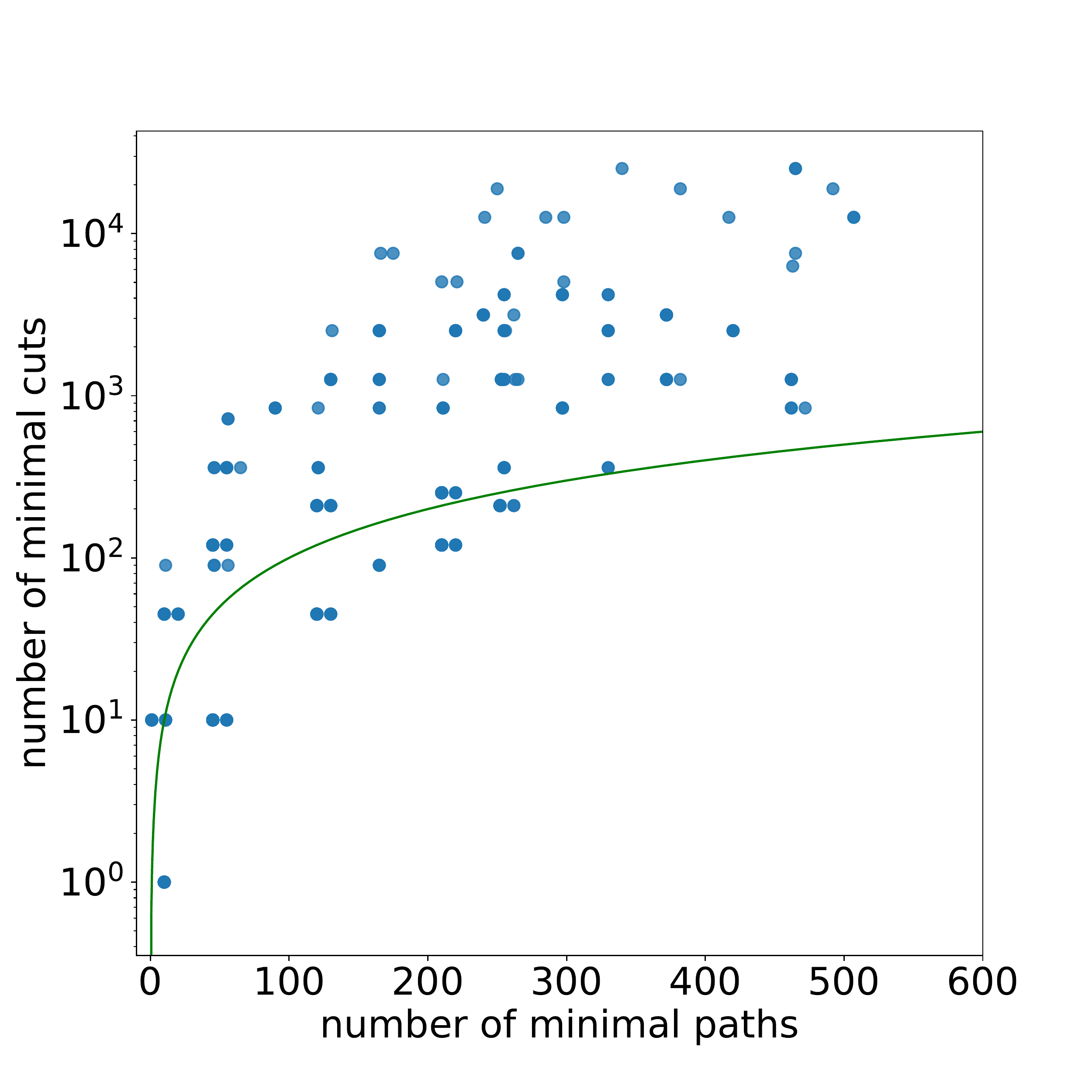} &   \includegraphics[scale=0.5]{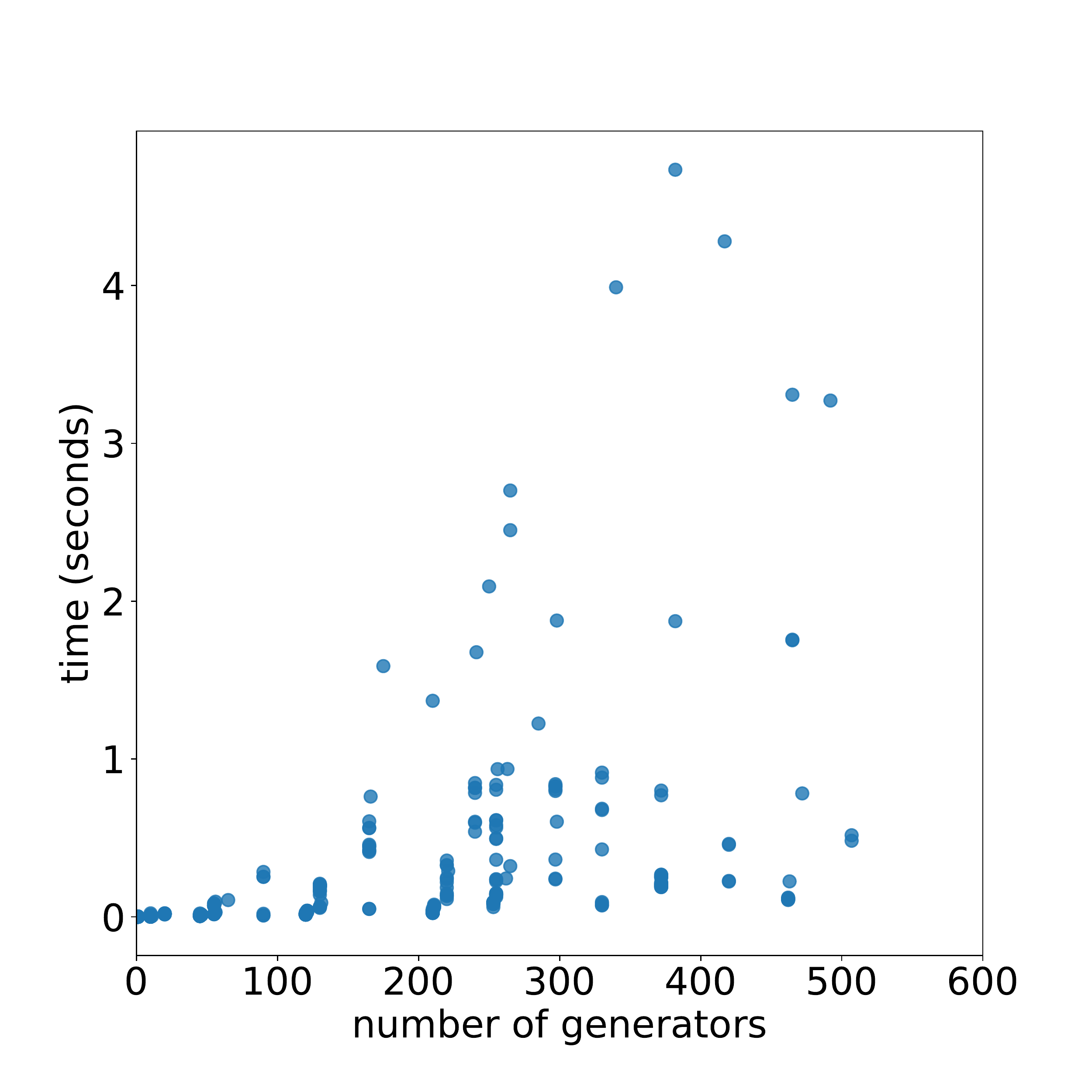}
 \\
(a)&(b) \\[6pt]
\end{tabular}}
\caption{Number of minimal paths vs. number of minimal cuts and computing time for generalized multi-state $k$-out-of-$n$ systems with $10$ components, $4$ levels, and non-identical probabilities.}
\label{fig:kn}
\end{figure}

\begin{figure}[b]
     \resizebox{\textwidth}{!}{
\begin{tabular}{cc}
 \includegraphics[scale=0.5]{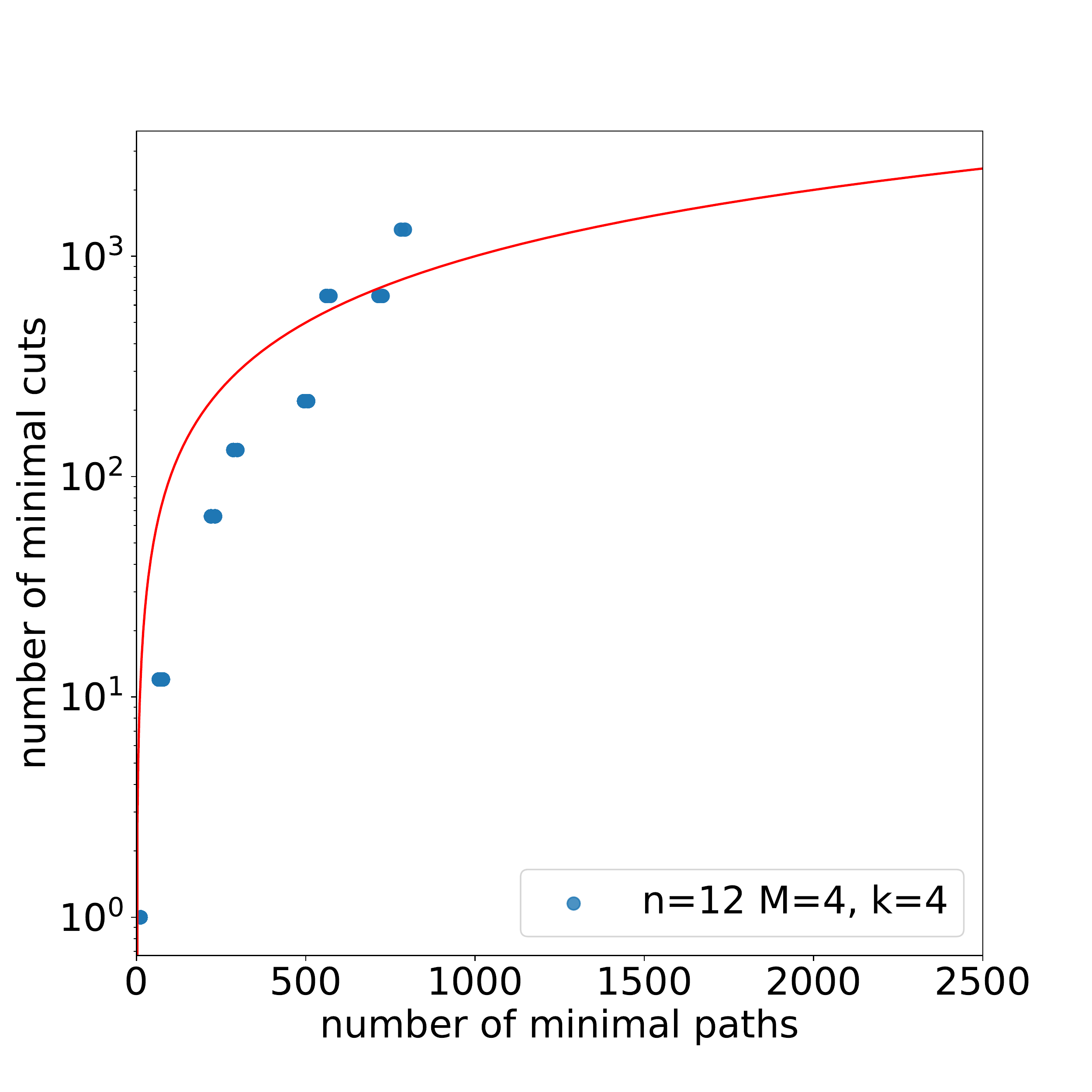} &   \includegraphics[scale=0.5]{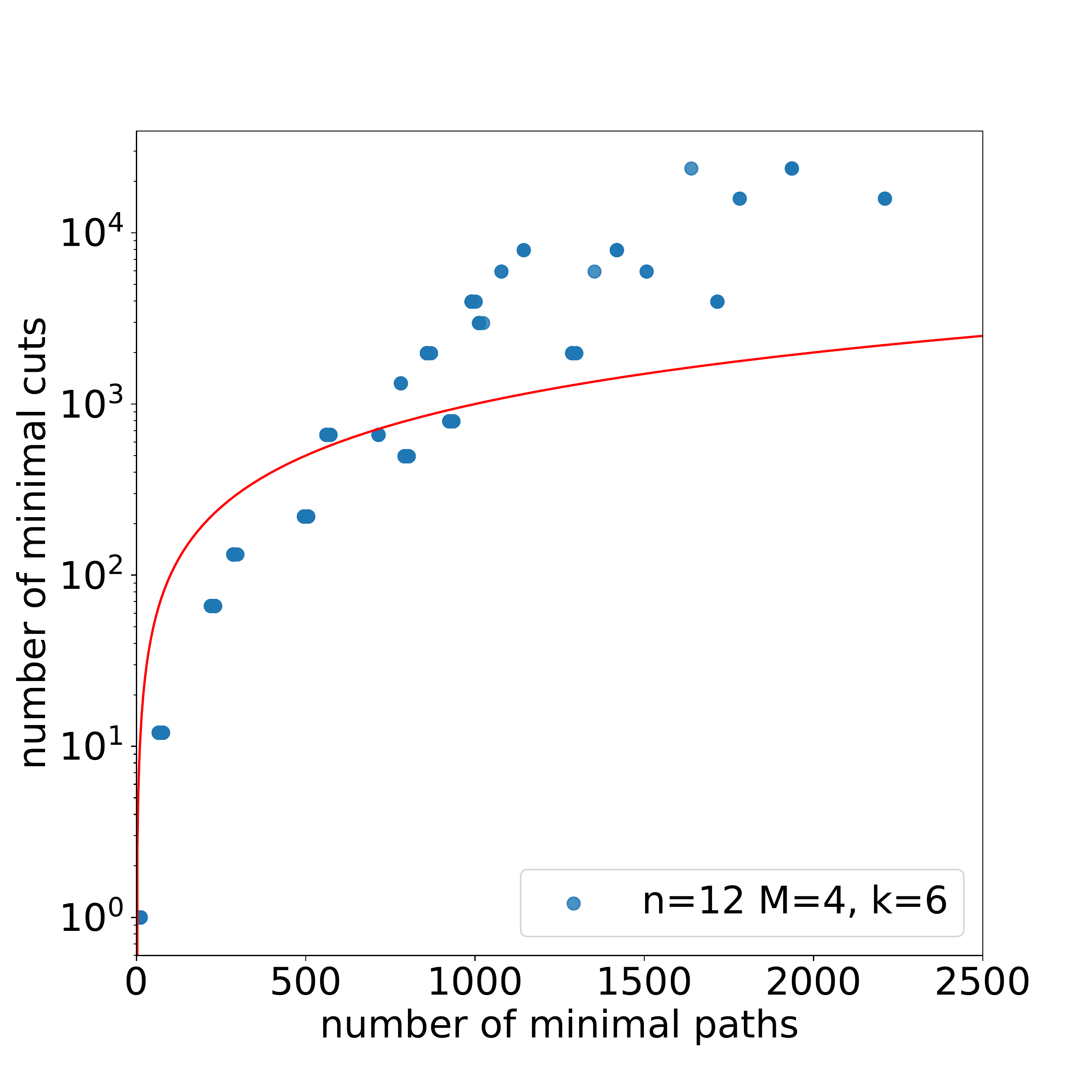}
 \\
  (a) &(b) \\[6pt]
  \multicolumn{2}{c}{\includegraphics[scale=0.5]{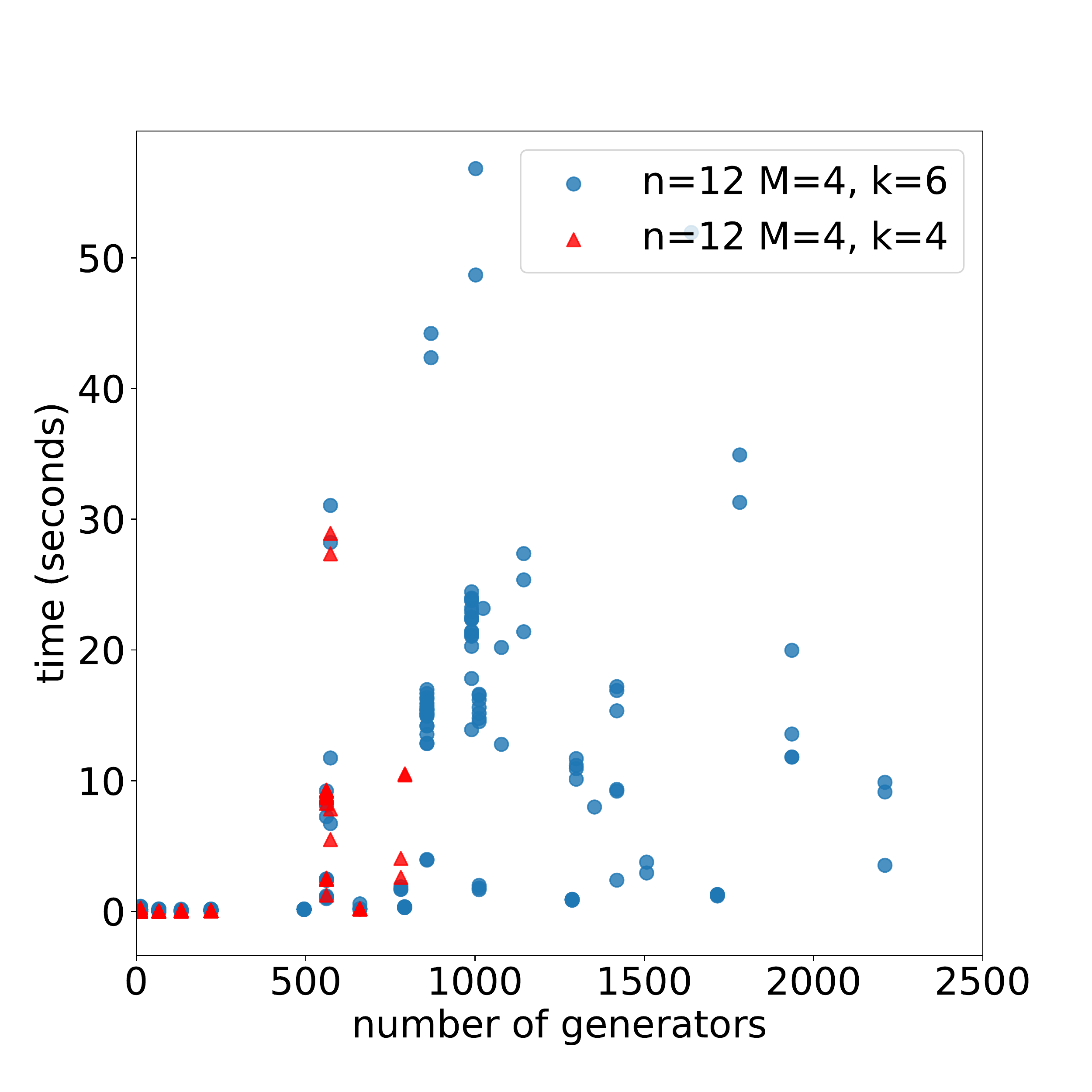}}\\
  \multicolumn{2}{c}{(c)}
\end{tabular}}
\caption{Number of minimal paths vs. number of minimal cuts and computing time for generalized multi-state $k$-out-of-$n$ systems with $12$ components, $4$ levels, $k=4,6$, and non-identical probabilities.}
\label{fig:v12l4}
\end{figure}

We have used our \verb|C++| class to compute the reliability of several generalized $k$-out-of-$n$ systems. Since each of these systems is given by a vector $(k_1,\dots,k_M)$ we generated randomly $100$ vectors for systems with four levels, and $10$ variables. Figures \ref{fig:kn} (a) and (b) show the number of minimal paths and minimal cuts of these systems and the computing time of these examples vs. the number of generators used for its computation in each case. All systems considered have components with independent, non-identical working probabilities. The figures show that most of these systems have a smaller number of minimal paths compared to the number of its minimal cuts, and that the computation time depends greatly on the structure of the system. Let us denote by $k$ the maximum of the integers $k_l$ for $l\in\{1,\dots,M\}$. Figures \ref{fig:v12l4} (a) and (b) show the number of minimal paths and cuts for systems with $12$ components, $4$ levels and $k=4$, $k=6$. The number of minimal cuts and paths of multi-state generalized $k$-out-of-$n$ grows as $n\choose k$. The performance of our algorithms depend greatly on the number of minimal paths or minimal cuts, as can be seen in Figure \ref{fig:v12l4} (c). There exist specialized algorithms for this kind of systems that are recursive on $M$, see \cite{ZT06,PSW20} or based on Decision Diagrams \cite{MLAD15}.

\subsection{Computational complexity}
The algebraic method is (in its general form) an enumerative method, similar to the inclusion-exclusion approach but less redundant. The compact form of the Hilbert series provided by our algorithms gives some computational advantages, but there exist certain intrinsic limitations due to the complexity of the problem. The computation of network reliability (either $k$-terminal, $2$-terminal or all-terminal) is a $\#P$-hard problem \cite{B86} and hence there is no hope of finding an efficient algorithm for computing the reliability of general systems unless $P=NP$, even in the binary case.
The algebraic method in which our algorithms are based shows that the problem of computing the reliability of a multi-state system can be polynomially reduced to the computation of the multigraded Hilbert series of a monomial ideal. This problem belongs to the class of $\#P$-hard problems and there exist several sub-problems of it that are known to be $\#P$-complete or $NP$-complete. In particular, the problem of computing the Euler characteristic of an abstract simplicial complex is equivalent to the computation of the coefficient of the monomial $x_1\cdots x_n$ in the multigraded Hilbert series of a (square-free) monomial ideal, and this problem belongs to the $\#P$-complete complexity class \cite{RS13}.

There are two complementary directions to follow for finding satisfactory solutions for these problems. One is to develop specialized polynomial algorithms for particular families of systems. The other is to find algorithms showing good heuristic behavior when applied to general problems. In these two directions it is of paramount importance to develop good implementations in terms of data types, memory management, etc. that make the algorithms applicable to practical problems.

The algebraic method for system reliability contributes to both of the directions described above. On the one hand, the study of the structure of the ideals of particular classes of systems can provide efficient algorithms or even formulas (explicit or recursive) for their Hilbert series, see \cite{SW10,SW11} for $k$-out-of-$n$ and consecutive $k$-out-of-$n$ binary systems, and \cite{PSW20} for the multi-state version. As an example, the formulas for $k$-out-of-$n$ systems have complexity $O(n^2)$ which is quadratic, but not optimal when restricted to systems with statistically independent components, for which the algorithm in \cite{B95}, based on the Fast Fourier Transform, runs in complexity $O(n(\log n)^2)$. On the other hand, for the general case we used efficient algorithms for computing the multi-graded Hilbert series of monomial ideals and Alexander duals. These algorithms avoid much of the redundancy that shows up in reliability computation of general systems, when we have no evident structure to take advantage of. Besides, they make use of the recursive nature of the problem, which has also been used in other approaches like the Universal Generating Function method. However, there is still room for improvement. As the UGF and other methods show, it is important, for the sake of efficiency, to identify good base cases for the recursion, and for simplification techniques. The algorithms provided in this paper use only {\em algebraic} base cases and simplifications, and hence it is expected to gain efficiency by exploring other base cases that arise from the knowledge of system reliability. This is beyond the scope of this paper and is left as future work. 
Finally, as it is common in computer algebra, implementations of general algorithms which are good enough for $NP$-hard or $\#P$-hard problems offer good performance in practice. A paradigmatic example of this are the good algorithms for Gr\"obner bases, a problem whose complexity is known to be doubly exponential. This is the case of the class presented in this paper, in which we took advantage of the data types and optimized routines provided by \verb|CoCoALib| together with good implementations for the Hilbert series and Alexander dual algorithms. This allows us to efficiently compute the reliability of general systems of big size with an affordable use of time and memory resources.

\section{Conclusions and further work}
We have presented a \verb|C++| class that computes the reliability of multi-state coherent systems. The class is included and distributed with the computer algebra library \verb|CoCoALib| as free software. The algorithms in this class are based on the algebraic approach to system reliability analysis developed in the last decade. The main advantage of these algorithms is that they can be applied to general systems, that they provide bounds and that can be applied without modification to systems with i.i.d probabilities and systems in which the probabilities are not identical. The main drawback is that this approach is enumerative, in the sense that relies on the enumeration of minimal paths or cuts, which may be impractical for big systems.

Specialized algorithms for particular systems are not easy to find but are very efficient in practice, see for example the algorithms based on Multi-Valued Decision Diagrams for multi-state $k$-out-of-$n$ systems \cite{MLAD15} or the linear algorithm for networks with small treewidth \cite{GM20}. Our future work includes the design of specialized algebraic algorithms for particular kinds of systems. The structure of particular systems induces a particular structure in the associated ideals which can be studied using algebraic and combinatorial tools allowing the design of more efficient algorithms, as described for instance in \cite{SW15,MPSW19}. Another direction of improvement is to adapt the algorithm for systems with non-independent components. The algebraic theory is exactly the same and only the probability assignment to the computed monomials need to be changed. This would give a wider flexibility to our \verb|C++| class. Finally, further tuning and optimization of the existing algorithms will likely improve their efficiency and reduce the computing times, in particular optimizations coming from the comparison and strong points of other methods, like the UGF method.

The fact that our approach is general makes it useful as one of the default algorithms to try in  the cases for which no specific algorithms are known yet, and also as a tool to benchmark new specific algorithms for such problems.

\section{Acknowledgements}
The first author was partially supported by the ``National Group for Algebraic and Geometric Structures, and Their Applications'' (GNSAGA, INdAM). The last two authors have been partially supported by grant MTM2017-88804-P from Ministerio de Econom\'ia, Industria y Competitividad (Spain). 

\bibliographystyle{ieeetr}
\bibliography{Classbibliography.bib}

\begin{thebibliography}{10}

\bibitem{reliasoft}
``Reliasoft.'' Available at \texttt{http://www.reliasoft.com}.

\bibitem{ALD}
``{ALD}.'' Available at \texttt{http://reliability-analysis-software.com/}.

\bibitem{ITEM}
``{ITEM} software.'' Available at \texttt{http://www.itemsoft.com}.

\bibitem{LC16}
H.-S. Li and Z.-J. Cao, ``Matlab codes of subset simulation for reliability
  analysis and structural optimization,'' {\em Struct. Multidisc. Optim.},
  vol.~54, pp.~391--410, 2016.

\bibitem{R15}
S.~Rao, {\em Reliability engineering}.
\newblock Pearson, 2015.

\bibitem{R20}
M.~Reid, ``Reliability. a {P}ython library for reliability engineering.''
  Available at \texttt{http://github.com/MatthewReid854/reliability}, 2020.

\bibitem{SYGBL08}
L.~B. Shaffer, T.~M. Young, F.~M. Guess, H.~Bensmail, and R.~V. Le\'on, ``Using
  {R} software for reliability data analysis,'' {\em International Journal of
  Reliability and Applications}, vol.~9, pp.~53--70, 2008.

\bibitem{sharpe}
R.~Sahner, K.~S. Trivedi, and A.~Pullafito, ``{SHARPE}, ({S}ymbolic
  {H}ierarchical {A}utomated {R}eliability and {P}erformance {E}valuator).''
  Available at \texttt{http://sharpe.pratt.duke.edu}.

\bibitem{tiobe}
``{TIOBE}.'' Available at \texttt{http://www.tiobe.com/tiobe-index}, april
  2020.

\bibitem{CoCoALib}
J.~Abbott and A.~M. Bigatti, ``{CoCoALib}: a {C}++ library for doing
  {C}omputations in {C}ommutative {A}lgebra.'' Available at
  \texttt{http://cocoa.dima.unige.it/cocoalib}.

\bibitem{BP88}
M.~Ball and J.~Provan, ``Computing network reliability in time polynomial in
  the number of cuts,'' {\em Oper. Res.}, vol.~32, pp.~516--526, 1988.

\bibitem{BP88b}
M.~Ball and J.~Provan, ``Disjoint products and efficient computations of
  reliability,'' {\em Oper. Res.}, vol.~36, pp.~703--715, 1988.

\bibitem{S91}
D.~Shier, {\em Network reliability and algebraic structures}.
\newblock Clarendon Press, 1991.

\bibitem{HKCD95}
D.~D. Harms, M.~Kratzel, C.~J. Colbourn, and J.~S. Devitt, {\em Network
  reliability. Experiments with a symbolic algebra environment}.
\newblock CRC Press, 1995.

\bibitem{U87}
I.~Ushakov, ``Optimal standby problem and a universal generating function,''
  {\em Sov. J. Comput. Sys. Sci.}, vol.~25, pp.~61--73, 1987.

\bibitem{L05}
G.~Levitin, {\em The Universal Generating Function in Reliability Analysis and
  Optimization}.
\newblock Springer, 2005.

\bibitem{LSX19}
D.~Lu, S.and~Shi and H.~Xiao, ``Reliability of sliding window systems with two
  failure modes,'' {\em Reliability Engineering \& System Safety}, vol.~188,
  pp.~366--376, 2019.

\bibitem{N11}
B.~Natvig, {\em Multi-state systems reliability theory with applications}.
\newblock John Wiley \& sons, 2011.

\bibitem{EPS78}
E.~El-Neweihi, F.~Proschan, and J.~Sethurman, ``Multi-state coherent systems,''
  {\em J. Applied Probability}, vol.~15, pp.~675--688, 1978.

\bibitem{KZ03}
W.~Kuo and M.~Zuo, {\em Optimal reliability modelling: principles and
  applications}.
\newblock John Wiley \& sons, 2003.

\bibitem{TB17}
K.~Trivedi and A.~Bobbio, {\em Reliability and availability engineering}.
\newblock Cambridge University Press, 2017.

\bibitem{GNW02}
B.~Giglio, D.~Q. Naiman, and H.~P. Wynn, ``Gr\"obner bases, abstract tubes, and
  inclusion–exclusion reliability bounds,'' {\em IEEE Trans. Rel.}, vol.~51,
  pp.~358--366, 2002.

\bibitem{GW04}
B.~Giglio and H.~P. Wynn, ``Monomial ideals and the scarf complex for coherent
  systems in reliability theory,'' {\em Annals of Statistics}, vol.~32,
  pp.~1289--1311, 2004.

\bibitem{D03}
K.~Dohmen, {\em Improved Bonferroni inequalities via abstract tubes}.
\newblock Springer, 2003.

\bibitem{SW09}
E.~{S\'aenz-de-Cabez\'on} and H.~P. Wynn, ``Betti numbers and minimal free
  resolutions for multi-state system reliability bounds,'' {\em Journal of
  Symbolic Computation}, vol.~44, pp.~1311--1325, 2009.

\bibitem{SW15}
E.~{S\'aenz-de-Cabez\'on} and H.~P. Wynn, ``Hilbert functions for design in
  reliability,'' {\em IEEE Trans. Rel.}, vol.~64, pp.~83--93, 2015.

\bibitem{MPSW19}
F.~Mohammadi, P.~Pascual-Ortigosa, E.~{S\'aenz-de-Cabez\'on}, and H.~Wynn,
  ``Polarization and depolarization of monomial ideals with application to
  multi-state system reliability,'' {\em Journal of Algebraic Combinatorics},
  vol.~51, pp.~617--639, 2020.

\bibitem{BK94}
R.~A. Boedigheimer and K.~C. Kapur, ``Customer-driven reliability models for
  multistate coherent systems,'' {\em IEEE Trans. Rel.}, vol.~43, pp.~46--50,
  1994.

\bibitem{PSW20}
P.~Pascual-Ortigosa, E.~{S\'aenz-de-Cabez\'on}, and H.~Wynn, ``Algebraic
  reliability of multi-state $k$-out-of-$n$ systems,'' {\em Probability in the
  Engineering and Informational Sciences}, 2020.

\bibitem{E95}
D.~Eisenbud, {\em Commutative algebra with a view towards algebraic geometry}.
\newblock Springer, 1995.

\bibitem{B97}
A.~M. Bigatti, ``Computation of {H}ilbert-{P}oincar\'e series,'' {\em Journal
  of Pure and Applied Algebra}, vol.~119, pp.~237--253, 1997.

\bibitem{MLAD15}
Y.~Mo, X.~Liudong, A.~S. V., and D.~J. B, ``Efficient analysis of multi-state
  k-out-of-n systems,'' {\em Reliability Engineering \& System Safety},
  vol.~133, pp.~95--105, 2015.

\bibitem{S09}
E.~{S\'aenz-de-Cabez\'on}, ``Multigraded betti numbers without computing
  minimal free resolutions,'' {\em Appl. Alg. Eng. Commun. Comput.}, vol.~20,
  pp.~481--495, 2009.

\bibitem{SC83}
A.~Satyanarayana and M.~K. Chang, ``Network reliability and the factoring
  theorem,'' {\em Networks}, vol.~13, pp.~107 -- 120, 1983.

\bibitem{MS04}
E.~Miller and B.~Sturmfels, {\em Combinatorial Commutative Algebra}.
\newblock Springer, 2004.

\bibitem{BS09}
A.~M. Bigatti and {S\'aenz-de-Cabez\'on}, ``Computation of the (n-1)-st
  {K}oszul {H}omology of monomial ideals and related algorithms,'' in {\em
  International Symposium on Symbolic and Algebraic Computation, ISSAC'09,
  Seoul, South Korea, 2009}, pp.~31--37, {ACM}, 2009.

\bibitem{CoCoA}
J.~Abbott, A.~M. Bigatti, and L.~Robbiano, ``{CoCoA}: a system for doing
  {C}omputations in {C}ommutative {A}lgebra.'' Available at
  \texttt{http://cocoa.dima.unige.it}.

\bibitem{Frobby}
B.~H. Roune, ``Frobby.'' Available at \texttt{http://www.broune.com/frobby}.

\bibitem{GN17}
J.~G\r{a}semyr and B.~Natvig, ``Improved availability bounds for binary and
  multistate monotone systems with independent component processes,'' {\em
  Journal of Applied Probability}, vol.~54, pp.~750--762, 2017.

\bibitem{FN85}
E.~Funnemark and B.~Natvig, ``Bounds for the availabilities in a fixed time
  interval for multistate monotone systems,'' {\em Adv. Appl. Prob}, vol.~17,
  pp.~638--665, 1985.

\bibitem{LL03}
A.~Lisnianski and G.~Levitin, {\em Multi-state system reliability}.
\newblock World Scientific, 2003.

\bibitem{LFD10}
A.~Lisnianski, I.~Frenkel, and Y.~Ding, {\em Multi-state system reliability
  analysis and optimization for engineers and industrial managers}.
\newblock Springer, 2010.

\bibitem{SXD10}
A.~Shrestha, L.~Xing, and Y.~Dai, ``Decision diagram based methods and
  complexity analysis for multi-state systems,'' {\em IEEE Trans. Rel.},
  vol.~59, pp.~145--160, 2010.

\bibitem{B79}
M.~Ball, ``Computing network reliability,'' {\em Oper. Res.}, vol.~27,
  pp.~823--838, 1979.

\bibitem{GM20}
A.~K. Goharshady and F.~Mohammadi, ``An efficient algorithm for computing
  network reliability in small trewidth,'' {\em Reliability Engineering \&
  System Safety}, vol.~193, 2020.

\bibitem{BTZ18}
G.~Bai, Z.~Tian, and M.~Zuo, ``Reliability evaluation of multistate networks:
  An improved algorithm using state space decomposition and experimental
  comparison,'' {\em IISE Transactions}, vol.~50, pp.~407--418, 2018.

\bibitem{ZLB20}
C.~Zhang, T.~Liu, and G.~Bai, ``An improved algorithm for reliability bounds of
  multistate networks,'' {\em Communications in Statistics - Theory and
  Methods}, 2020.

\bibitem{ER59}
P.~Erd\H{o}s and A.~R\'enyi, ``On random graphs i.,'' {\em Publicationes
  Mathematicae Debrecen}, vol.~6, pp.~290--297, 1959.

\bibitem{BA99}
A.-L. Barab\'asi and R.~Albert, ``Emergence of scaling in random networks,''
  {\em Science}, vol.~286, pp.~509--512, 1999.

\bibitem{HZW00}
J.~Huang, M.~J. Zuo, and Y.~Wu, ``Generalized multi-state k-out-of-n:g
  systems,'' {\em IEEE Trans. Rel.}, vol.~49, pp.~105--111, 2000.

\bibitem{ZT06}
M.~J. Zuo and Z.~Tian, ``Performance evaluation of generalized multi-state
  k-out-of-n systems,'' {\em IEEE Trans. Rel.}, vol.~55, pp.~319--327, 2006.

\bibitem{B86}
M.~Ball, ``Computational complexity of network reliability analysis: An
  overview,'' {\em IEEE Trans. Rel.}, vol.~R-35, pp.~230--239, 1986.

\bibitem{RS13}
B.~H. Roune and E.~{S\'aenz-de-Cabez\'on}, ``Complexity and algorithms for the
  euler characteristic of simplicial complexes,'' {\em Journal of Symbolic
  Computation}, vol.~50, pp.~170--196, 2013.

\bibitem{SW10}
E.~{S\'aenz-de-Cabez\'on} and H.~P. Wynn, ``Betti numbers and minimal free
  resolutions for multi-state system reliability bounds,'' {\em Appl. Alg. Eng.
  Commun. Comput.}, vol.~21, pp.~443--457, 2010.

\bibitem{SW11}
E.~{S\'aenz-de-Cabez\'on} and H.~P. Wynn, ``Computational algebraic algorithms
  for the reliability of generalized $k$-out-of-$n$ and related systems,'' {\em
  Math. Comput. Simulation}, vol.~82, pp.~68--78, 2011.

\bibitem{B95}
L.~Belfore, ``An $o(n(\log n)^2)$ algorithm for computing the reliability of
  $k$-out-of-$n:g$ and \& $k$-to$l$-out-of-$n:g$ systems,'' {\em IEEE Trans.
  Rel.}, vol.~R-44, pp.~132--136, 1995.

\end{thebibliography}

\newpage
\section*{Appendix: UML diagram of the Algebraic Reliability C++ class}
\begin{figure}[h]
  \centering
  \includegraphics[scale=0.50]{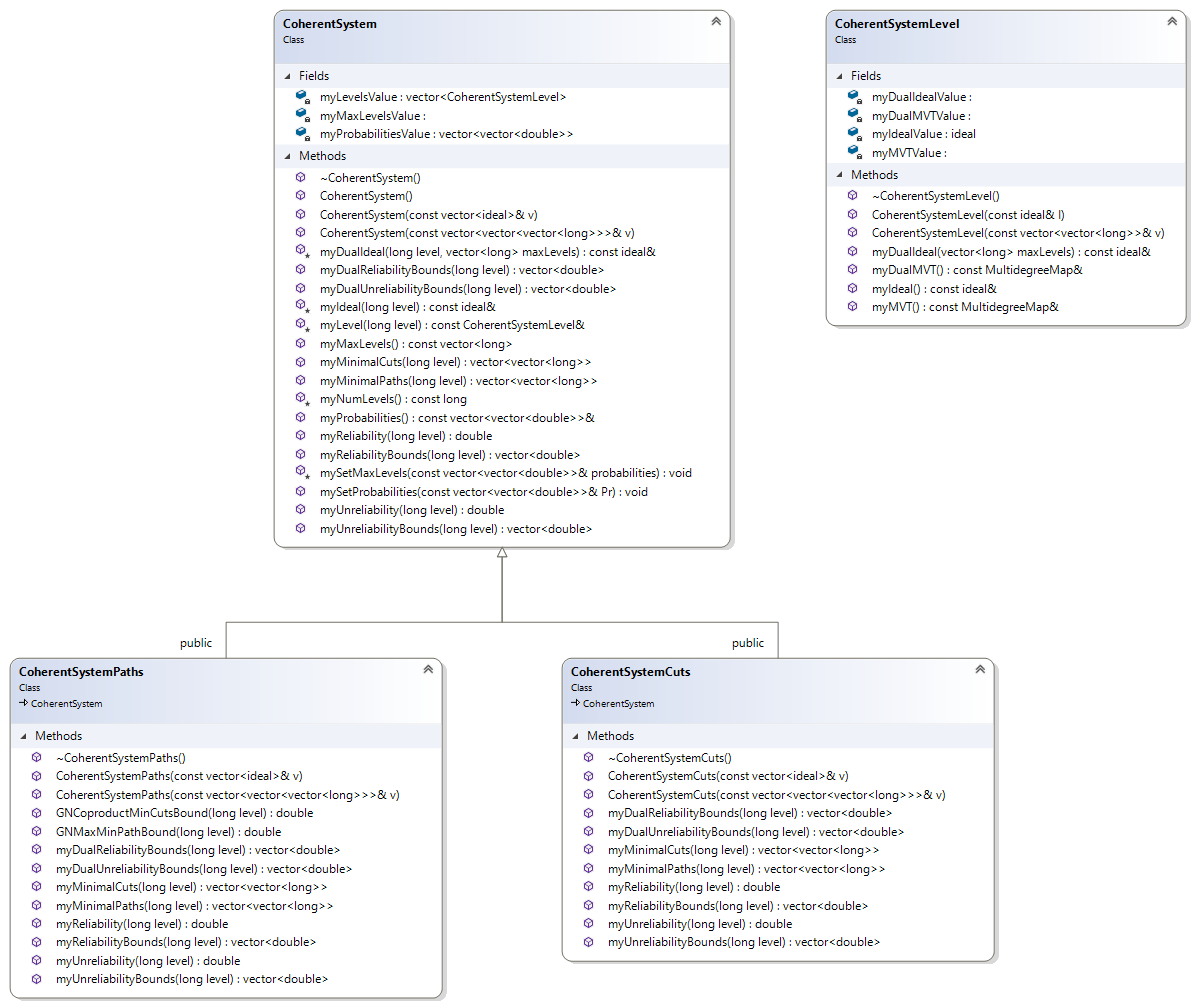}
  \caption{\texttt{UML} diagram of the \texttt{CoherentSystem} class}\label{fig:classDiagram}
\end{figure}

\end{document}